\documentclass[11pt,reqno,a4paper]{amsart}
\usepackage{color}
\usepackage{amssymb,amsmath}
\usepackage[latin1]{inputenc}
\usepackage[active]{srcltx}
\usepackage{graphicx}
\usepackage{exscale,relsize}
\usepackage{textgreek}
\usepackage{epsfig,graphics}
\usepackage{psfrag}
\usepackage{caption}
\usepackage{subcaption}
\usepackage[toc,page]{appendix}
\usepackage{bigints}
\def\N{\mathbb{N}}
\def\R{\mathbb{R}}
\def\m1{{I\!\!M}}

\renewcommand{\to}{\rightarrow}
\newcommand{\pa}{\partial}

\newcommand{\ino}{\int_{\Omega}}

\newcommand{\ainf}{\mbox{as\;}\;n\to+\infty}

\newcommand{\rife}[1]{(\ref{#1})}
\newcommand{\ov}[1]{\overline{#1}}
\newcommand{\un}[1]{\underline{#1}}

\newcommand{\scp}{\scriptstyle}
\newcommand{\sscp}{\scriptscriptstyle}

\renewcommand{\dfrac}{\displaystyle\frac}
\newcommand{\finedim}{\hspace{\fill}$\square$}
\newcommand{\intbar}{\mathop{\int\makebox(-15.5,0){\rule[6pt]{.7em}{0.3pt}}\kern-6pt}\nolimits}

\newcommand{\ii}{\infty}

\newcommand{\eps}{\varepsilon}
\newcommand{\dt}{\delta}

\newcommand{\al}{\alpha}
\newcommand{\bvl}{{\boldsymbol \eta}_{\sscp \lm}}
\newcommand{\vlv[1]}{\eta_{#1,\sscp \lm}}
\newcommand{\alv[1]}{\alpha_{#1,\sscp \lm}}
\newcommand{\plv[1]}{\psi_{#1,\sscp \lm}}
\newcommand{\ulv[1]}{u_{#1,\sscp \lm}}

\newcommand{\tlv[1]}{\tau_{#1,\sscp \lm}}

\newcommand{\rlav[1]}{\mbox{\Large \textrho}_{#1,\sscp \lm,\sscp \al_{#1}}}
\newcommand{\rlavv[1]}{V_{#1,\sscp \lm,\sscp \al_{#1}}}

\newcommand{\rlv[1]}{\mbox{\Large \textrho}_{#1,\sscp \lm}}
\newcommand{\rlqv[1]}{V_{#1,\sscp \lm}}
\newcommand{\mlv[1]}{m_{#1,\ssl}}
\newcommand{\bal}{{\boldsymbol\alpha}_{\sscp \lm}}

\newcommand{\bmu}{{\boldsymbol\mu}}
\newcommand{\ba}{{\boldsymbol\alpha}}
\newcommand{\bvarp}{{\boldsymbol\varphi}}
\newcommand{\bF}{{\mathbf F}}
\newcommand{\bs}{{\mathbf s}}
\newcommand{\bphi}{{\boldsymbol\phi}}
\newcommand{\bpl}{{\boldsymbol\psi}_{\sscp \lm}}

\newcommand{\bp}{{\boldsymbol\psi}}

\newcommand{\sg}{\sigma}

\newcommand{\om}{\Omega}
\newcommand{\lm}{\lambda}

\newcommand{\ur}{\underline{\rho}}
\newcommand{\urls}{\underline{\rho}_{\sscp \lm}}
\newcommand{\urlu}{\underline{\rho}_{1}}
\newcommand{\urld}{\underline{\rho}_{2}}

\newcommand{\rh}{\mbox{\Large \textrho}}

\newcommand{\pl}{\psi_{\sscp \lm}}

\newcommand{\ul}{u_{\sscp \lm}}

\newcommand{\ssl}{\sscp \lm}

\newcommand{\all}{\al_{\ssl}}

\newcommand{\els}{E_{\ssl,s}}

\newcommand{\el}{E_{\ssl}}

\newcommand{\fbi}{{\bf (F)$_{\lm}$}}
\newcommand{\prl}{{\textbf{(}\mathbf P\textbf{)}_{\mathbf \lm}}}

\newtheorem{theorem}{Theorem}[section]
\newtheorem{proposition}[theorem]{Proposition}
\newtheorem{lemma}[theorem]{Lemma}
\newtheorem{corollary}[theorem]{Corollary}
\newtheorem{remark}[theorem]{Remark}
\newtheorem{definition}[theorem]{Definition}
\newcommand{\brm}{\begin{remark}\rm}
\newcommand{\erm}{\end{remark}}
\newcommand{\bdf}{\begin{definition}\rm}
\newcommand{\edf}{\end{definition}}
\newcommand{\bte}{\begin{theorem}}
\newcommand{\ete}{\end{theorem}}
\newcommand{\bpr}{\begin{proposition}}
\newcommand{\epr}{\end{proposition}}
\newcommand{\ble}{\begin{lemma}}
\newcommand{\ele}{\end{lemma}}
\newcommand{\bco}{\begin{corollary}}
\newcommand{\eco}{\end{corollary}}
\newcommand{\beq}{\begin{equation}}
\newcommand{\eeq}{\end{equation}}
\newcommand{\bdm}{\begin{displaymath}}
\newcommand{\edm}{\end{displaymath}}

\newcommand{\graf}[1]{\left\{\begin{array}{ll}#1\end{array}\right.}

\def\sideremark#1{\ifvmode\leavevmode\fi\vadjust{\vbox to0pt{\vss
 \hbox to 0pt{\hskip\hsize\hskip1em \vbox{\hsize2.1cm\tiny\raggedright\pretolerance10000 \noindent #1\hfill}\hss}\vbox to15pt{\vfil}\vss}}}

\begin{document}
\numberwithin{equation}{section}
\parindent=0pt
\hfuzz=2pt
\frenchspacing

\title[Uniqueness for Lane-Emden free boundary type systems]{
A Lane-Emden system of free boundary type:\\ existence, uniqueness and monotonicity of solutions}


\author[D. Bartolucci]{Daniele Bartolucci}
\address{Daniele Bartolucci, Department of Mathematics, University of Rome \emph{"Tor Vergata"}, Via della ricerca scientifica n.1, 00133 Roma.}
\email{bartoluc@mat.uniroma2.it}

\author[Y. Hu]{Yeyao Hu}
\address{Yeyao Hu, School of Mathematics and Statistics, HNP-LAMA, Central South University, Changsha,
Hunan 410083, P. R. China.}
\email{huyeyao@gmail.com}

\author[A. Jevnikar]{Aleks Jevnikar}
\address{Aleks Jevnikar, Department of Mathematics, Computer Science and Physics, University of Udine, Via delle Scienze 206, 33100 Udine, Italy.}
\email{aleks.jevnikar@uniud.it}

\author[J. Wei]{Juncheng Wei}
\address{Juncheng Wei, Department of Mathematics, Chinese University of Hong Kong, Shatin, Hong Kong.}
\email{wei@math.cuhk.edu.hk}

\author[W. Yang]{Wen Yang}
\address{Wen Yang, Department of Mathematics, Faculty of Science and Technology, University of Macau, Macau, P.R. China}
\email{wenyang@um.edu.mo}

\thanks{D.B. is partially supported by
MIUR Excellence Department Project MatMod@TOV awarded to the Department of Mathematics, Univ. of Rome Tor Vergata, by PRIN 2022
"\emph{Variational and Analytical aspects of Geometric PDEs}" and by E.P.G.P. Project funded by Univ. of Rome Tor Vergata 2024. A.J. is partially supported by PRIN 2022 "\emph{Pattern formation in nonlinear phenomena}". D.B. and A.J. are members of the INDAM Research Group GNAMPA. Y.H. is supported by NSFC China No. 12471205. J.W. is partially supported by GRF fund of RGC of Hong Kong entitled "New frontiers in singularity formations of nonlinear partial differential equations." W.Y. is supported by the National
Key R\&D Program of China 2022YFA1006800, NSFC China No.\,12171456, NSFC China No.\,12271369, FDCT No. 0070/2024/RIA1, Start-up
Research Grant No.\,SRG2023-00067-FST,
Multi-Year Research Grant No.\,MYRG-GRG2024-00082-FST and UMDF-TISF/2025/006/FST. }

\begin{abstract}
We consider a Hamiltonian system of free boundary type, showing first uniform bounds and existence of solutions and of the free boundary. Then, for any smooth and bounded domain, we prove uniqueness
of positive solutions in a suitable interval and
show that the associated energies and boundary values have a monotonic behavior.
Some consequences are discussed about the parametrization of the unbounded Rabinowitz continuum for a
class of superlinear strongly coupled elliptic systems.
\end{abstract}
\maketitle

\tableofcontents


\setcounter{section}{0}
\setcounter{equation}{0}
\section{\bf Introduction}

Let $\om\subset \R^N$, $N\geq2$, be an open and bounded domain of class $C^{3}$,
we are concerned with the following constrained Hamiltonian system of free boundary type,
$$
\graf{-\Delta v_1 = \lm (v_2)_{+}^{p_2}\quad \mbox{in}\;\;\om\\ \\
-\Delta v_2 = \lm (v_1)_{+}^{p_1}\quad \mbox{in}\;\;\om\\ \\
-\bigintss\limits_{\pa\om} \frac{\pa v_1}{\pa\nu}=1=-\bigintss\limits_{\pa\om} \frac{\pa v_2}{\pa\nu} \\ \\
v_1=\al_1,\quad v_2=\al_2 \quad \mbox{on}\;\;\pa\om
}\qquad \qquad \mbox{\bf (F)}_{\lm}
$$
for the unknowns $\al_i \in \R$, and
$v_i\in C^{2,r}(\ov{\om}\,)$, $i=1,2$, for some fixed $r\in (0,1)$. Here, $(v)_+$ is the positive part of $v$, $\nu$ is the exterior unit normal, $\lm> 0$ and
\beq\label{Souto}
\frac{1}{p_1+1}+\frac{1}{p_2+1}>\frac{N-2}{N-1},\quad p_i\in (0,+\ii),\,i=1,2.
\eeq
The relevance of the limiting hyperbola in \eqref{Souto} for classical Hamiltonian elliptic systems of Lane-Emden type was first noticed in \cite{Souto}, see \cite{defig} for more details.
Putting,
$$
p_{N}=\graf{+\ii,\;N=2 \vspace{0.2cm}\\ \frac{N}{N-2}\,,\; N\geq 3,}
$$
we remark that, as far as $N\geq3$,  $p_1=p_{N}=p_2$ satisfy $\frac{1}{p_1+1}+\frac{1}{p_2+1}=\frac{N-2}{N-1}$.

To simplify the exposition, by a suitable scaling of $\lm$ we
assume that $|\om|=1$.\\

The system \fbi\, is a vectorial generalization of the classical (\cite{Te,Te2,BeBr}) "scalar" free boundary problem which is obtained
in the particular case $p_1=p_2$, $\al_1=\al_2$, $v_1=v_2$,
whose study is motivated by Tokamak's plasma physics (\cite{Fre,Wes}). In the scalar case
$p_{N}$ turns out to be a natural critical exponent see \cite{BeBr,Ort}.\\
Another motivation to purse the analysis
of \fbi\, is to find a parametrization of solutions of the Hamiltonian strongly
coupled elliptic system,
$$
\graf{-\Delta u_1 = \mu_2 (1+u_2)^{p_2}\quad \mbox{in}\;\;\om\\
-\Delta u_2 = \mu_1 (1+u_1)^{p_1}\quad \mbox{in}\;\;\om\\
u_1>0,\quad u_2>0 \quad \mbox{in}\;\;\om\\
u_1=0,\quad u_2=0 \quad \mbox{on}\;\;\pa\om.
}\qquad \qquad \mbox{\bf (H)}
$$
For $\mu_1=\mu_2$ the existence of an unbounded continuum of "scalar" ($u_1=u_2$) solutions of ${\bf (H)}$ follows from the classical result in \cite{Rab}.
The analysis of Hamiltonian
elliptic systems is a classical subject and we refer to \cite{defig} for a comprehensive introduction about this topic, see also \cite{Soup} and references
therein.
Minimal solutions branches (in the sense of Crandall-Rabinowitz (\cite{CrRab})) and
multiplicity results for general systems including ${\bf (H)}$ has been described in \cite{CC1,CC2}, see also
\cite{mont} and references therein.\\
However our main motivation comes from the fact that
we are not aware of any result either about uniqueness and qualitative behavior of branches of solutions
of \fbi\, or just about the qualitative behavior of \underline{non minimal} solutions of ${\bf (H)}$.
Remark that if $p_1\neq p_2$ then \fbi\, has no scalar solutions. In particular concerning ${\bf (H)}$,
inspired by recent results in \cite{BHJY}, \cite{BJ}, \cite{BJ1}, we look for integral quantities naturally arising from $\prl$ to describe
the monotonic behavior of the solutions.\\

First of all let us consider the auxiliary problem,
$$
\graf{-\Delta \psi_1 =(\al_2+{\lm}\psi_2)_+^{p_2}\quad \mbox{in}\;\;\om\\ \\
-\Delta \psi_2 =(\al_1+{\lm}\psi_1)_+^{p_1}\quad \mbox{in}\;\;\om\\ \\
\bigintss\limits_{\om} (\al_2+{\lm}\psi_2)_+^{p_2}=1=\bigintss\limits_{\om} (\al_1+{\lm}\psi_1)_+^{p_1}\\ \\
\psi_i=0 \quad \mbox{on}\;\;\pa\om,\quad i=1,2, \\ \\
\al_i\in\R,\quad i=1,2,
}\qquad \prl
$$
for the unknowns $\al_i\in\R$ and $\psi_i \in C^{2,r}_{0}(\ov{\om}\,)$, $i=1,2$. Here $\lm\geq 0$,
$(p_1,p_2)$ satisfy \eqref{Souto} and, for some fixed $r\in (0,1)$, we set,
$$
C^{2,r}_0(\ov{\om}\,)=\{\psi \in C^{2,r}(\ov{\om}\,)\,:\, \psi=0\mbox{ on }\pa \om\},\;
C^{2,r}_{0,+}(\ov{\om}\,)=\{\psi \in C^{2,r}_0(\ov{\om}\,)\,:\, \psi> 0\mbox{ in } \om\}.
$$

Here and in the rest of this paper we refer to solutions of $\prl$ of this sort as \un{classical solutions}.\\
Interestingly enough, problem $\prl$ seems to be of independent interest as it defines the stationary solutions
in the study of two species chemiotaxis models with nonlinear diffusion recently pursued in \cite{CKL}, see also
Remark \ref{rem1.1} below.

For $\lm\geq 0$ fixed it is useful to denote a solution of $\prl$ by
$\bal=(\alv[1],\alv[2])$, $\bpl=(\plv[1],\plv[2])$, and define positive/non negative solutions as follows,\\

{\bf Definition.}
{\it We say that $(\bal,\bpl)$ is a positive/non negative solution of {\rm $\prl$} if $\alv[i]>0$/$\alv[i]\geq 0$,
$i=1,2$, respectively.}\\

Clearly, if $(\bal,\bpl)$ is a positive solution, then by the strong maximum principle
$\plv[i]>0$ in $\om$, $i=1,2$. Remark that for $\lm>0$, $(\bal,\bpl)$ is a
solution of {\rm $\prl$} if and only if $((\alv[1],\alv[2]),(v_{1,\sscp \lm},v_{2,\sscp \lm}))$
solves \fbi,\,with $v_{i,\sscp \lm}=\alv[i]+\lm \plv[i]$, $i=1,2$.\\

\bigskip

We point out that, since $|\om|=1$ and $\lm\geq 0$ by assumption, then if $(\bal,\bpl)$ is a non negative
solution of $\prl$ then necessarily,
$$
\alv[i]\leq 1,\;i=1,2,
$$
and the equality $\alv[1]=1=\alv[2]$ holds if and only if $\lm=0$. We will frequently use this fact without further comments.
Actually, if $\lm=0$, then $\prl$ takes the form
$$
\graf{-\Delta \psi_{1, 0} =1\quad \mbox{in}\;\;\om\\ \\
-\Delta \psi_{2, 0} =1\quad \mbox{in}\;\;\om\\ \\
\psi_{1, 0}=0=\psi_{1, 0} \quad \mbox{on}\;\;\pa\om
}
$$
and admits a unique (in fact scalar) solution
$(\ba_0, \bp_0) =((1,1),(G[1],G[1]))$,
where we define,
$$
G[\rh](x)=\ino G_{\om}(x,y)\rh(y)\,dy,\;x\in\om.
$$
Here $G_{\om}$ is the Green function of $-\Delta$ with Dirichlet boundary conditions on $\om$.
Obviously, to say that $(\bal,\bpl)$ is a solution of $\prl$ is the same as to say that
$\bpl=(G[\rlv[2]],G[\rlv[1]])$ and $\ino \rlv[1]=1=\ino \rlv[2]$,
where, unless otherwise specified, we set
$$
\rlv[i]=(\al_{i,\lm}+\lm\psi_{i,\lm})_+^{p_i},\quad i=1,2.
$$

\medskip

By a standard fixed point argument (see Appendix \ref{appF}) it can be shown that, for any $\lm>0$ small enough,
there exists at least one solution of $\prl$ and in particular that $\al_i>\frac13$, $i=1,2$ for any such
solution of $\prl$.  Also, by a well known argument in \cite{BeBr}, one could prove the existence of at least one solution
for any $\lm>0$ as far as $p_i<p_N$, $i=1,2$.\\
By using the weak Young inequality (\cite{Lieb}), we refine here the variational argument in \cite{BeBr}, see section \ref{appD},
to come up with at least one solution of $\prl$ whenever $(p_1,p_2)$ satisfy \eqref{Souto}. Remark that, still
as far as \eqref{Souto} is satisfied, we can prove that if $(\psi_{1,\sscp \lm},\psi_{2,\sscp \lm})\in W^{2,p_2}_0(\om)\times W^{2,p_1}_0(\om)$ is just assumed to be a strong solution of $\prl$, then $(\psi_{1,\sscp \lm},\psi_{2,\sscp \lm})\in C^{2,r_0}_0(\om)\times C^{2,r_0}_0(\om)$ and satisfy
to some uniform bound for bounded $\lm$, see Lemma \ref{lemE1} below. These uniform estimates seems to be new and in particular are crucial for our purposes. Also, at least in the scalar case $p_1=p=p_2$, $\psi_1=\psi_2$,
they are sharp since in fact, if $p\geq p_{N}$, it is well known (\cite{WYe},\cite{We}) that solutions may blow up for $\lm$ large enough.
It is also not too difficult to prove that our variational functional (see \eqref{jeil} below) is in fact not anymore
coercive as far as $p_1=p=p_2$ and $p\geq p_{N}$.
\brm\label{rem1.1} {\it After the completion of this work, we came to learn about the recent reference \cite{CKL} where essentially the same variational
functional is analyzed {\rm (see \eqref{jeil} below)} on the whole space $\R^N$\!. It seems to be an interesting open problem to extend our uniform estimates {\rm (Lemma \ref{lemE1} below)} to the larger region defined as follows,
$$p_1\left(p_2-\frac{2}{N-2}\right)<p_{N}\quad \mbox{or}\quad p_2\left(p_1-\frac{2}{N-2}\right)<p_{N},$$ as far as $N\geq 3$.
Although not explicitly used in \cite{CKL}, these inequalities follow just by considering the range of parameters pursued therein. Also, we do not exclude that some arguments in \cite{CKL} could be used to come up with the existence of
solutions of {\rm $\prl$} in this larger region. Remark that the intersection point of these two hyperbolas is the symmetric boundary point $(p_{N},p_{N})$ in \eqref{Souto}.\\
 More in general, it could be interesting to investigate the relevance of the well known critical hyperbolas pushed forward in
\cite{CDM1} and \cite{CDM2} for problem {\rm $\prl$}.}
\erm

\bigskip

However, our main concern is about uniqueness and qualitative behavior of $(\bal,\bpl)$
depending on $\lm$. This is not trivial for two reasons. First of all we have in principle four unknowns to control,
which are $\al_{i,\lm}$ and $\psi_{i,\lm}$, $i=1,2.$ On the other side, due to the constraints in $\prl$,
as recently proved in (\cite{BJ1}) in the scalar case $\pl=\plv[1]=\plv[2]$, $\all=\alv[1]=\alv[2]$,
for $\lm>0$ small enough $\all$ is strictly decreasing while, by standard
arguments (\cite{CrRab}), for \un{fixed} $\all=\al$ and disregarding the constraints in $\prl$,
then $\pl$ is strictly increasing for any $\lm$ small enough. Actually the same monotonicity property holds,
of course at fixed $(\al_1,\al_2)$ and for $\lm$ small enough, for $\plv[1]$ and $\plv[2]$, due to well known results
about the maximum principle for cooperative elliptic systems (\cite{dFM}). Therefore, unlike classical scalar problems (\cite{CrRab}) there is
a competition between the monotonic behavior of $(\al_{i,\lm}+\lm \psi_{i,\lm})$ as a function of $\lm$.\\
This is why we do not adopt classical maximum principles based argument but rather rely on ideas recently pursued in \cite{BJ1,BJW},
see also \cite{B2,BJ,BW} where different class of problems are considered.

We will prove existence, uniqueness and monotonicity via a
refined dual spectral formulation suitable to analyze positive solutions of the constrained problem $\prl$.
The first eigenvalue in this
spectral setting is denoted by $\sg_1(\bal,\bpl)$, see section \ref{sec2}. In particular, we will prove the
monotonicity
of two naturally defined variational
quantities associated to $\prl$ (see section \ref{appD}), which are the energy,
$$
\el:=\ino\rlv[1] G[\rlv[2]]\equiv \ino\rlv[i] \plv[i]=\ino(\nabla \plv[1],\nabla \plv[2]),\,i=1,2,
$$
and the free energy,
$$
F_{\ssl}=\frac{1}{r_1}\ino (\rlv[1])^{r_1}+\frac{1}{r_2}\ino (\rlv[2])^{r_2}-\lm
\ino \rlv[1] G[\rlv[2]],
$$
where $r_i=1+\frac{1}{p_i}$, $i=1,2$. Furthermore, we can prove the monotonicity of the linear combination
$$
\dfrac{p_1\alv[1]}{p_1+1}+  \dfrac{p_2\alv[2]}{p_2+1}\,.
$$

Set ${\mathbf p}=(p_1,p_2)$ and
$$
\lm^*(\om,{\mathbf p})=\sup\{\lm>0\,:\,\sg_1(\bal,\bpl)>0,\,\al_{i,\mu}>0,\,i=1,2,\mbox{\rm \,for any solution of }
{\textbf{(}\mathbf P\textbf{)}_{\mathbf \mu}},\,\forall\,\mu<\lm\}.
$$
It can be shown, see Lemma \ref{lmsmall} in Appendix \ref{appF} and Proposition \ref{prsigma}
in section \ref{sec2}, that $\lm^*(\om,{\mathbf p})$
is well defined and strictly positive and our first task is to prove that $\lm^*(\om,{\mathbf p})<+\ii$.
More exactly our first result is about the existence of a free boundary
in the interior of $\om$ for solutions of \fbi\, with $(p_1,p_2)$ satisfying \eqref{Souto}.
The point here is that one would like to know whether or not, for a fixed $\lm$, $\min\{\alv[1],\alv[2]\}$
is negative, which implies in particular that at least one among
$\om_{i,-}:=\{x\in\om\,:\,v_{i}<0\}$, $i=1,2$ is not empty.
In the scalar case this problem for $p=1$ is fully understood,
see \cite{BeBr,Pudam,Te2}, while, for $p>1$, the existence of a multiply connected free boundary
has been proved in \cite{We} for $\lm$ large
and under some assumptions about the existence of non degenerate critical points of a suitably defined
{Kirchoff-Routh} type functional. Still for $\lm$ large, but only for $N=2$ and for domains with
non trivial topology, a similar result has been obtained in \cite{Liu}.
Other sufficient conditions for the existence of solutions with $\al<0$ has been found in \cite{AmbM},
which however assume the nonlinearity $v_+^p$ to be replaced by $g_+(x,v)$ satisfying $g(x,t)\geq ct$,
for some $c>0$, which therefore does not fit our scalar problem. More recently it has been shown in
\cite{BHJY} that there are no scalar
solutions with $\all\geq 0$ for $\lm$ large. We generalize that result here to the case of the
cooperative and strongly coupled systems of free boundary type \fbi.

\bte\label{thlambda} Let $(p_1,p_2)$ satisfy \eqref{Souto}. Then we have:

\begin{itemize}
\item[$(a)$] \emph{(Existence)} For any $\lm>0$ there exists at least one solution of \fbi.

\medskip

\item[$(b)$] \emph{(Existence of free boundary)} Suppose either $N=2$ or $N\geq3$ and $\Omega$ convex. Then, there exists
$\ov{\lm}=\ov{\lm}(\om,{\mathbf p},N)>0$ depending only on ${\mathbf p}, N$ and $\om$ such that if
$(\bal,\bpl)$ is a non negative
solution of {\rm$\prl$}, then $\lm\leq \ov{\lm}$. In particular, for any $\lm>\ov{\lm}$ we have $\min\{\alv[1],\alv[2]\}<0$ for any solution of \fbi.
\end{itemize}
\ete

As mentioned above the existence part in $(a)$ follows by a refinement of the variational argument for the scalar case about $\prl$ provided in
\cite{BeBr}. We prove $(b)$ by a blow up argument, which is a well known tool in the study of a priori estimates
for Lane-Emden systems, see \cite{defig} and references therein. However the situation here is slightly different from standard
models, which is why we provide a self contained proof.\\

We denote by $\mathcal{G}(\om)$ the set of solutions of
$\prl$ for $\lm\in [0,\lm^*(\om,{\mathbf p}))$. Our next aim is to show uniqueness of these solutions, see Theorem~\ref{thmLE} below. Observe that, under the assumption of Theorem~\ref{thlambda}-${(b)}$, we have $\lm^*(\om,{\mathbf p})<+\ii$. Actually, since we do not expect uniqueness for any $\lambda\in(0,+\ii)$, we believe Theorem \ref{thlambda}-(b) holds true also for non-convex domains with $N\geq3$.

Let $B_r=\{x\in \R^N\,:\, |x|<r\}$ with volume $|B_r|$,
$\mathbb{D}_{\sscp N}$ be the $N$-dimensional ball of unit volume and let us denote by,
\beq\label{sob}
\Lambda(\om,t)=\inf\limits_{w\in H^1_0(\om), w\equiv \!\!\!\!/ \;0}\dfrac{\ino |\nabla w|^2}{\left(\ino |w|^{t}\right)^{\frac2t}}\,,
\eeq
which provides the best constant in the Sobolev embedding $\|w\|_p\leq S_p(\om)\|\nabla w\|_2$,
$S_p(\om)=\Lambda^{-2}(\om,p)$, $p\in[1,2p_{N})$. Let us set
$$
\sg_{1,*}=\liminf\limits_{\lm\to \lm^*(\om,{\mathbf p})^-}\sg_1(\bal,\bpl),
\quad
\al_{i,*}=\liminf\limits_{\lm\to \lm^*(\om,{\mathbf p})^-}\alv[i],\,i=1,2.$$

For the sake of clarity, we point out that if a map $\mathcal{M}$ from an interval $[a,b]\subset \R$ to a
Banach space $X$ is said to be real analytic,
then it is understood that $\mathcal{M}$ can be extended in
an open neighborhood of $a$ and $b$ where it admits a power
series expansion, totally convergent in the $X$-norm.\\
Here $\|\cdot\|_{i,\ssl}$ stands for the weighted norms naturally associated to the problem, see section \ref{sec2} for definitions.

\bte\label{thmLE}
Let $(p_1,p_2)$ satisfy \eqref{Souto}. Then we have:
\begin{itemize}
\item[1.] \emph{(Uniqueness)} For any $\lm \in [0,\lm^*(\om,{\mathbf p}))$
there exists a unique solution $(\bal,\bpl)$ of {\rm $\prl$} and
$\mathcal{G}(\om)$ is a real analytic simple curve of positive solutions
$[0,\lm^*(\om,{\mathbf p}))\ni \lm\mapsto (\bal,\bpl)$. As $\lm\to 0^+$ we have, \\
$
\bal=(1,1)+\mbox{\rm O}(\lm),\;\;\bpl=(\psi_{1,0},\psi_{2,0})+\mbox{\rm O}(\lm),
\el=E_0(\om)+\mbox{\rm O}(\lm),
$
where,
$$
E_{0}(\om)=\ino \ino G_{\om}(x,y)\,dxdy\leq
E_{0}(\mathbb{D}_{\sscp N})=
{\textstyle\dfrac{|B_1|^{\scp -\frac 2 N}}{2(N+2)}}.
$$

\item[2.] \emph{(Monotonicity)} For any
$\lm\in [0,\lm^*(\om,{\mathbf p}))$ it holds,
\beq\label{monot}
\frac{d F_{\ssl}}{d\lm}<0,\;\frac{d \el}{d\lm}\geq 0,\;\frac{d}{d\lm}
\left(\frac{p_1\alv[1]}{p_1+1}+\frac{p_2\alv[2]}{p_2+1}\right)<0.
\eeq
Moreover,
$$
\frac{d \el}{d \lm} \geq p_1\|[\plv[1]]_{1,\ssl}\|^2_{1,\ssl}+p_2\|[\plv[2]]_{2,\ssl}\|^2_{2,\ssl}.
$$

\item[3.] \emph{(Spectral estimates)} If either $\sg_{1,*}=0$ or if $\al_{i,*}=0$, $i=1,2$, then
$$\lm^*(\om,{\mathbf p})\geq \frac{1}{p_2}\Lambda(\om,2p_2).$$
\end{itemize}
\ete

\bigskip

Clearly $E_0(\om)$ is just the torsional rigidity of $\om$.
The above theorem holds for any smooth and bounded domain,
in any dimension and for any subcritical (in the sense of \eqref{Souto}) exponent. Remark that the result is sharp
in the scalar case $p=p_1=p_2$, $\all=\alv[1]=\alv[2]$,  $\pl=\plv[1]=\plv[2]$,
where we have that $\all> 0$ and $\sg_1(\all,\pl)>0$ for any
$\lm<\frac{1}{p}\Lambda(\om,2p)$, see \cite{BJ1}. In particular, it has been shown
in \cite{BJ1} that if $p<p_N$ there exists a positive solution for any
$\lm< \frac{1}{p}\Lambda(\om,2p)$,
then by Theorem \ref{thmLE} we immediately deduce the following corollary about the case $1\leq p_1= p_2<p_{N}$.

\bco Let $1\leq p_1= p_2<p_{N}$. For any $\lm \in [0,\lm^*(\om,{\mathbf p}))$ the unique solution of
{\rm $\prl$} is the scalar solution $\all=\alv[1]=\alv[2]$,  $\pl=\plv[1]=\plv[2]$.
\eco

\bigskip

We still do not know whether or not $\sg_{1,*}=0$.
However, as far as $p_1\neq p_2$, it seems that
$\al_{i,*}=0$, $i=1,2$ is not a natural assumption, actually a more reasonable guess is that, in general,
if $\sg_{1,*}>0$,
then either $\al_{1,*}=0,\al_{2,*}>0$ or
$\al_{1,*}>0,\al_{2,*}=0$. This is interesting since in this case we come up with the parametrization of
an unbounded branch of solutions of ${\bf (H)}$.\\
For any classical solution ${\bf u}=(u_1,u_2)$ of ${\bf (H)}$ for some $\bmu=(\mu_1,\mu_2)\in ([0,+\ii))^2$ we define,

$$
\gamma(\bmu,{\bf u})=\frac{p_1}{p_1+1}\frac{1}{\|1+u_1\|_{p_1}}+\frac{p_2}{p_2+1}\frac{1}{\|1+u_2\|_{p_2}},
$$

$$
E(\bmu,{\bf u})=\frac{1}{2\mu_1\mu_2}
\ino\frac{\mu_1(1+u_1)^{p_1}u_1+\mu_2(1+u_2)^{p_2}u_2}{\|1+u_1\|^{p_1}_{p_1}\|1+u_2\|^{p_2}_{p_2} },
$$

$$
F(\bmu,{\bf u})=\gamma(\bmu,{\bf u})+\frac{p_1p_2-1}{(p_2+1)(p_1+1)}
E(\bmu,{\bf u}).
$$

Then we have,

\bte\label{thmH} Let $(p_1,p_2)$ satisfy \eqref{Souto} and $\mathcal{G}(\om)$ be the set of unique
solutions $(\bal,\bpl)$ of {\rm $\prl$} for $\lm \in [0,\lm^*(\om,{\mathbf p}))$.
Then
$$
{\bf u}_{\sscp \lm}=(\ulv[1],\ulv[2])=\left(\frac{\lm}{\alv[1]}\,\plv[1],\frac{\lm}{\alv[2]}\,\plv[2]\right),
$$
is a solution of ${\bf (H)}$ with,
$$
\bmu_{\sscp\lm}=(\mu_{1,\ssl},\mu_{2,\ssl})=
\left(\lm\,\frac{\alv[1]^{p_1}}{\alv[2]},\lm\,\frac{\alv[2]^{p_2}}{\alv[1]}\right),
$$
and for any
$\lm\in [0,\lm^*(\om,{\mathbf p}))$ it holds,
\beq\label{monot-u}
\frac{d F}{d\lm}(\bmu_{\ssl}, {\bf u}_{\ssl})<0,\quad \frac{d E}{d\lm}(\bmu_{\ssl}, {\bf u}_{\ssl})\geq 0,
\quad\frac{d \gamma}{d\lm}
(\bmu_{\ssl}, {\bf u}_{\ssl})<0.
\eeq

In particular, if $\al_{1,*}=0,\al_{2,*}>0$ then
$(\bmu_{\sscp\lm}, {\bf u}_{\sscp \lm})$
is a real analytic and unbounded curve and,  possibly along a subsequence, we have that
$$
\mu_{1,\ssl}\to 0, \quad \mu_{2,\ssl}\to +\ii, \quad \|1+u_{1,\ssl}\|_{p_1}\to +\ii.
$$
Moreover, under the assumption of Theorem \ref{thlambda}-$(b)$, we have
$\frac{u_{1,\ssl}}{\|1+u_{1,\ssl}\|_{p_1}}\to \lm^*\psi_{1,*}, \frac{u_{2,\ssl}}{\|1+u_{2,\ssl}\|_{p_2}}
\to \lm^*\psi_{2,*}$, with convergence in $C^2(\ov{\om})$
where $(\psi_{1,*},\psi_{2,*})$ is a solution of {\rm $\prl$} with $\lm=\lm^*=\lm^*(\om,{\mathbf p})$,
$\al_1=0$ for some $\al_2>0$. The conclusion is analogous in the case $\al_{1,*}>0,\al_{2,*}=0$.
\ete

\bigskip

Concerning Theorem \ref{thmLE}, due to the competition between the monotonic behavior of $\alv[i]$ and $\plv[i]$,
it seems difficult to attack the problem by arguments based on the standard maximum principle. On the other side,
the fact that $\el$, $F_{\ssl}$ are monotonic increasing
could be deduced for variational solutions of $\prl$ once we know the uniqueness
of solutions. However this is not enough to claim the smoothness of $\el$ or the monotonicity/smoothness of
$\left(\frac{p_1\alv[1]}{p_1+1}+\frac{p_2\alv[2]}{p_2+1}\right)$.
The problem is more subtle for $\el$ since it seems that variational arguments do not yield any
information in this case.\\

Therefore, the proof of Theorem \ref{thmLE} relies on the interplay between a refined
bifurcation analysis and the variational formulation of $\prl$.\\
The crucial point is the set up of a sort of dual
 "Hamiltonian" spectral theory for the linearized operator of $\prl$,
 see the definition \rife{eLl}
of ${\mathbf L}_{\ssl}=(L_{1,\ssl},L_{2,\ssl})$ and the related eigenvalues equation \rife{lineq0}
in section \ref{sec2}. Remark that the eigenvalue equation is a non standard one, which is why we
refer to it as an "Hamiltonian" eigenvalue problem. The use of the operator ${\mathbf L}_{\ssl}$
is rather delicate also because it arises as the linearization of a vectorial constrained problem
$(\ino \rlv[i]=1, i=1,2)$ with respect to $(\bal,\bpl)$, which yields a non-local problem. As a consequence
it is not true in general, neither for scalar solutions, that its first eigenvalue, which we denote by $\sg_1(\bal,\bpl)$,
is simple and neither that
if $\sg_1(\bal,\bpl)$ is positive then the maximum principle holds. For example this is exactly
what happens in the scalar case for $\lm=0$ on $\mathbb{D}_{\sscp 2}$, where $\sg_1(\al_{\sscp 0},\psi_{\sscp 0})$ can be
evaluated explicitly (see \cite{BJ}) and one finds that
$\sg_1(\al_{\sscp 0},\psi_{\sscp 0})=\lm^{(2,0)}(\mathbb{D}_{\sscp 2})\simeq \pi (3,83)^2$ has three eigenfunctions,
two of which indeed change sign. Here
$\lm^{(2,0)}(\om)$ is the first non vanishing eigenvalue of $-\Delta$ on $\om$ on
the space of $H^1(\om)$ vanishing mean functions
with constant boundary trace. See also \cite{courant} for a related results.\\

However, if for a positive solution $(\bal,\bpl)$ with $\lm\geq 0$ it holds
$0\notin \sigma({\mathbf L}_{\ssl})$,
where $\sigma({\mathbf L}_{\ssl})$ stands for the spectrum of ${\mathbf L}_{\ssl}$,
then by the real analytic implicit function theorem (\cite{but}) the set of solutions of $\prl$
is locally a real analytic curve of positive solutions. In particular a real analytic curve of
positive solutions exists around $(\ba_0, \bp_0)$. Since solutions of $\prl$ are uniformly bounded
(see Lemma \ref{lemE1}) then for any $\lm<\lm^*(\om,{\mathbf p})$
there exists a unique solution and these solutions form a real analytic curve which we
denoted by $\mathcal{G}(\om)$.
An a priori bound from below far away from zero
for $\sg_1(\bal,\bpl)$ for positive solutions can be derived at this stage,
see Proposition \ref{prsigma}. An estimate
about the range where the $\al_{i,\ssl}$, $i=1,2$ may possibly vanish at the same time follows as well,
see Proposition \ref{prsigma1}.
At this point, since we know that $F_{\ssl}$, $\el$ and $\alv[i]$ are real analytic as functions of $\lm$ as far as
$\lm<\lm^*(\om,{\mathbf p})$, then, by the variational characterization
of $(\bal,\bpl)$, we deduce the monotonicity properties of  $F_{\ssl}$, $\el$ and $\alv[i]$.
The estimate about the derivative of $\el$
requires a more careful analysis of the Fourier expansion of $\frac{d \el}{d\lm}$ in terms of the
"Hamiltonian" Fourier basis, see section \ref{appD}.

It would be interesting to find a fourth monotonic quantity naturally associated to the problem, for example the self-interaction energy
$$
\els:=\frac12 \ino\rlv[2] G[\rlv[2]] + \frac12 \ino\rlv[1] G[\rlv[1]].
$$
However, it seems not easy to catch the qualitative behavior of $\els$ and we postpone this problem to a future work.

\bigskip

This paper is organized as follows. In section \ref{sec1} we discuss about the existence of the free boundary and prove Theorem \ref{thlambda}.
In section \ref{sec2} we set up the spectral and bifurcation analysis with the needed spectral
estimates. In section \ref{appD} we prove existence of variational solutions, uniqueness and monotonicity,
which yield the proof of Theorem \ref{thmLE} and, as a corollary, that of Theorem \ref{thmH}.
The proof of the existence of solutions and of the positivity of $\al_{i,\ssl}$ for $\lm$ small is
discussed in Appendix \ref{appF}.

\bigskip
\bigskip

\section{{\bf Existence of the free boundary}}\label{sec1}
For later purposes, see either Theorem \ref{thm5.2} below, we prove a regularity result of independent interest, showing that if
$(\bal,\bpl)$ is a solution of $\prl$ such that $(\psi_{1,\sscp \lm},\psi_{2,\sscp \lm})\in W^{2,p_2}_0(\om)\times W^{2,p_1}_0(\om)$ is just assumed
to be a strong solution of $\prl$, then $(\psi_{1,\sscp \lm},\psi_{2,\sscp \lm})\in C^{2,r_0}_0(\om)\times C^{2,r_0}_0(\om)$ and is
uniformly bounded in $C^{2,r_0}_0(\om)\times C^{2,r_0}_0(\om)$, as far as $\lm$ is bounded as well.
Indeed we have,
\ble\label{lemE1} Let ${\mathbf p}=(p_1,p_2)$ satisfy \eqref{Souto}. For any $\ov{\lm}>0$ there
exists a positive constant $C_1=C_1(r,\om,\ov{\lm},{\bf p},N)$ depending only on $\om$, $\ov{\lm}$, ${\bf p}$, $N$
and $r\in (0,1)$ such that $\|\plv[i]\|_{C^{2,r_0}_0(\ov{\om})}\leq C_1$, $i=1,2$ for any strong solution
$(\bal,\bpl)$ of {\rm $\prl$} with $\lm\in [0,\ov{\lm}\,]$, where $r_0=\min\{p_1,p_2,r\}$.
\ele
\proof
To simplify the notations we set $(\al_i,\psi_i)=(\alv[i],\plv[i])$.
Since $\om$ is of class $C^3$, by standard elliptic estimates
and a bootstrap argument it is enough to prove that either $\psi_1$ or $\psi_2$ is bounded. We prove only the case $N\geq 3$ which is more delicate.\\

Let $(p_1,p_2)$ satisfy\eqref{Souto} and assume w.l.o.g. that $p_1<p_{N}$.
Since $\ino (\al_i+\lm\psi_i)_+^{p_i}=1$, $i=1,2$, then it is well known (\cite{St4}) that for any $t\in [1,\frac{N}{N-1})$ there exists $C=C(t,N,\om)$ such that $\|\psi_i\|_{W_0^{1,t}(\om)}\leq C(t,N,\om)$,
$i=1,2$ for any solution of $\prl$. Thus, by the Sobolev inequality, for any $1\leq s< \frac{N}{N-2}$ we have
$\|\psi_i\|_{L^s(\om)}\leq C(s,N,\om)$, $i=1,2$ for some $C(s,N,\om)$. By the maximum principle $\psi_i\geq 0$, $i=1,2$.
Thus, either $\al_i>0$ and then $\ino (\al_i+\lm \psi_i)=1$ and (recall $|\om|=1$) $\al_i\leq 1$, or $\al_i<0$ and then $(\al_i+\psi_i)_+\leq \psi_i$. Since $p_1<p_{N}$,
then for any $1<m<\frac{p_{N}}{p_1}$, $\|(\al_1+\lm\psi_1)^{p_1}\|_{L^m(\om)}\leq C(p_1,N,{\lm},s,\om)$.\\
From now on we suppress the indications of the properties of the various constants involved in the estimates, being understood that $C$ is just a suitable uniform constant which do not depend by the solutions.\\

By standard elliptic theory
we have that $\|\psi_2\|_{W_0^{2,m}(\om)}\leq C$, for any $m<\frac{p_N}{p_1}$ and by the Sobolev embedding either $p_1\leq \frac{2}{N-2}$ and then
$\|\psi_2\|_{L^{r}(\om)}\leq C$ for any $r\geq 1$ and then the desired conclusion follows in a standard way by a bootstrap argument,
or $\frac{2}{N-2}<p_1<p_{N}$ and then $\|\psi_2\|_{L^{r}(\om)}\leq C$ for any $r<\frac{N}{p_1-\frac{2}{N-2}}$. By \eqref{Souto} we have
$$
p_2<\frac{p_1+N}{p_1(N-2)-1},
$$
and observe that, putting $p_1=\frac{x}{N-2}, x\in (2,N)$,
$$
m_{N,p_{1}}:=\frac{\frac{N}{p_1-\frac{2}{N-2}}}{\frac{p_1+N}{p_1(N-2)-1}}=\frac{N(N-2)}{p_1+N}\frac{p_1(N-2)-1}{p_1(N-2)-{2}}=(N-2)f_{_N}(x),
$$
where $f_{_N}(x)=\frac{x-1}{(x-2)(1+a_{_N}x)},a_{_N}=\frac{1}{N(N-2)}$. Therefore, since $f_{_N}$ is decreasing, we see that $m_{N,p_{1}}$ is
monotonic decreasing, with $m_{N,p_{1}}\to +\ii$ as $p_1\to \left(\frac{2}{N-2}\right)^+$, $m_{N,p_{1}}\to (N-2)^+$ as $p_1\to \left(\frac{N}{N-2}\right)^-$.
In particular $\|(\al_2+\lm\psi_2)^{p_2}\|_{L^n(\om)}\leq C$ for any $n<m_{N,p_{1}}$, whence by standard elliptic theory
$\|\psi_1\|_{W_0^{2,n}(\om)}\leq C$, for any $n<m_{N,p_{1}}$ and either $m_{N,p_{1}}\geq \frac{N}{2}$ and then
$\|\psi_1\|_{L^{r}(\om)}\leq C $ for any $r\geq 1$ and then the desired conclusion follows in a standard way by a bootstrap argument,
or $m_{N,p_{1}}< \frac{N}{2}$ and then $\|\psi_1\|_{L^{r}(\om)}\leq C $ for any
$r<\frac{N}{N-2 m_{N,p_{1}}}$. On the other side, as deduced above, we have that $m_{N,p_{1}}\geq (N-2)$ whence in particular $m_{N,p_{1}}\geq \frac{N}{2}$, as far $N\geq 4$ and we are just left with the case $N=3$ which requires a different argument.\\
With the notations adopted above, we have that  $\|(\al_2+\lm\psi_2)^{p_2}\|_{L^n(\om)}\leq C$ for any $n<m^{(1)}_{3,p_{1}}$, where
$$
m^{(1)}_{3,p_{1}}=\frac{1}{1+\frac{p_1}{3}}\frac{p_1-1}{p_1-{2}} \quad \mbox{ and }\quad 1< m^{(1)}_{3,p_{1}}<\frac{3}{2}.
$$
Here the assumption about $p_1$ is that $p_1\in (\ov{p},3)$ where $m^{(1)}_{3,\ov{p}}=\frac{3}{2}$  (actually $\ov{p}> \frac{5}{2}$), where we recall that $m^{(1)}_{3,p_{1}}$ is decreasing in $p_1$ and $m_{3,3}=1$.
In particular, as mentioned above, by standard elliptic theory $\|\psi_1\|_{W^{2,n}_0(\om)}\leq C$ for any $n<m^{(1)}_{3,p_{1}}$ and by the Sobolev
embedding $\|\psi_1\|_{L^{s}(\om)}\leq C$ for any $s<s_2:=\frac{3}{3-2m^{(1)}_{3,p_{1}}}$. We define
$$
\sigma_1:=\frac{s_2}{3}=\frac{1}{3-2m^{(1)}_{3,p_{1}}}>1, \quad \delta:=\sg_1-1,
$$
whence $\|(\al_2+\lm\psi_1)^{p_1}\|_{L^m(\om)}\leq C$ for any $m<\frac{s_2}{p_1}=\frac{3}{p_1}\sg_1$.
Therefore, assuming w.l.o.g. that
$$
\sg_1<\frac{p_1}{2},
$$
by standard elliptic theory $\|\psi_2\|_{W^{2,m}_0(\om)}\leq C$ for any $m<\frac{3}{p_1}\sg_1$ and by the Sobolev embedding
$\|\psi_2\|_{L^{r}(\om)}\leq C$ for any $r<\frac{3\sg_1}{p_1-2\sg_1}$.
As a consequence
$\|(\al_2+\lm\psi_2)^{p_2}\|_{L^n(\om)}\leq C$ for any $n<m^{(2)}_{3,p_1}$ where, by using \eqref{Souto} once more, we define,
$$
m^{(2)}_{3,p_1}=\frac{3\sg_1}{p_1-2\sg_1}\frac{p_1-1}{p_1+3}=\sg_1\frac{p_1-2}{p_1-2\sg_1}m^{(1)}_{3,p_1}>\sg_1 m^{(1)}_{3,p_1}=\frac{m^{(1)}_{3,p_1}
}{3-2m^{(1)}_{3,p_{1}}}.
$$
Obviously we can assume w.l.o.g. that $m^{(2)}_{3,p_1}<\frac{3}{2}$, that is, using the last equality,
$$
1<m^{(1)}_{3,p_1}<\frac{9}{8} \;\mbox{ and consequently } 1<\sg_1<\frac{4}{3}.
$$
Therefore, as above
$\|\psi_1\|_{L^s(\om)}\leq C $ for any $s<s_3:=\frac{3}{3-2m^{(2)}_{3,p_{1}}}$ and we define
$$
\sigma_2:=\frac{s_3}{s_2}=\frac{3-2m^{(1)}_{3,p_{1}}}{3-2m^{(2)}_{3,p_{1}}}.
$$
At this point, since $m^{(2)}_{3,p_{1}}>\sg_1 m^{(1)}_{3,p_{1}}$, elementary arguments show that
$$
\sigma_2>\frac{3-2m^{(1)}_{3,p_{1}}}{3-2\sg_1m^{(1)}_{3,p_{1}}}=1+\frac{2m^{(1)}_{3,p_{1}}(\sg_1-1)}{3-2\sg_1m^{(1)}_{3,p_{1}}}>1+\dt=\sg_1.
$$
Consequently  $\|(\al_2+\lm\psi_1)^{p_1}\|_{L^m(\om)}\leq C$ for any $m<\frac{s_3}{p_1}=\frac{3}{p_1}\sg_1\sg_2>\frac{3}{p_1}\sg_1^2$ and by
standard elliptic theory and the Sobolev embedding we have that $\|(\al_2+\lm\psi_2)^{p_2}\|_{L^n(\om)}\leq C$ for any $n<m^{(3)}_{3,p_1}$ where, by using \eqref{Souto} once more, we define,
$$
m^{(3)}_{3,p_1}:=\frac{3\sg_1^2}{p_1-2\sg_1^2}\frac{p_1-1}{p_1+3}=\sg_1\frac{3\sg_1}{p_1-2\sg_1^2}\frac{p_1-1}{p_1+3}>\sg_1 m^{(2)}_{3,p_1},
$$
where we can assume w.l.o.g. that,
$$
\sg^2_1<\frac{p_1}{2}.
$$
At this point, by induction, it is not too difficult to prove that after a finite number of iterations, either  $p_1-2\sg_1^{k-1}\leq 0$ or $m^{(k)}_{3,p_1}=\frac{3\sg_1^{k-1}}{p_1-2\sg_1^{k-1}}\frac{p_1-1}{p_1+3}$ will become larger than $\frac{3}{2}$ and the desired conclusion follows, in this case as well.
\finedim

\bigskip
\bigskip

Next we present the proof of Theorem \ref{thlambda} about a priori estimates and the existence of a free boundary.\proof
We postpone the proof of $(a)$, i.e. the existence of at least one solution of $\prl$
for any $\lm>0$, to Theorem \ref{thm5.2} in section \ref{appD}.\\
Therefore we are going to prove $(b)$.
We argue by contradiction and assume that there exists
a sequence of non negative solutions $((\al_{1,n},\al_{2,n}),(\psi_{1,n},\psi_{2,n}))$ of $\left.\prl\right|_{\lm=\lm_n}$ such that
$\lm_n\to +\ii$. By Lemma \ref{lemE1}, for any fixed $n$ it holds $\|\psi_{i,n}\|_{\ii}\leq C_n$, $i=1,2$.
Let $m_{i,n}=\sup\limits_{\om}(\al_{i,n}+\lm_n\psi_{i,n})$, $i=1,2$. We split the proof in various steps.\\

STEP 1. Along a subsequence we have  $\max\{m_{1,n},m_{2,n}\}\to +\ii$, $\ainf$.\\
We will need the following lemma whose proof can be found in \cite{BHJY}.
\ble\label{minen}
Let $\psi$ be any solution of
$$
\graf{-\Delta \psi=f \quad \mbox{in}\;\;\om\\ \psi=0 \quad \mbox{on}\;\;\pa \om}
$$
where $\ino f=1$ and $\ino |f|^N\leq C$. Then,
$$
\ino|\nabla \psi|^2\geq c>0,
$$
for some positive constant $c>0$ depending only by $C$, $N$ and $\om$.
\ele

At this point we argue by contradiction and assume that there exists $C>0$ such that
$$
\sup\limits_n\max\{m_{1,n},m_{2,n}\}\leq C,
$$
so that
$\sup\limits_n\max\{\|\psi_{1,n}\|_\ii,\|\psi_{2,n}\|_\ii\}\leq \frac{C}{\lm_n}$. Therefore, along a subsequence
we have,

$$
\ino |\nabla \psi_{1,n}|^2=\ino (\al_{2,n}+\lm_n\psi_{2,n})^{p_2}\psi_{1,n}\leq \frac{C}{\lm_n}\to 0,\,\ainf,
$$

which contradicts Lemma \ref{minen} since $f_{i,n}=(\al_{i,n}+\lm_n\psi_{i,n})^{p_i}$
obviously satisfies the needed assumptions for any $i=1,2$.\\
Therefore, along a subsequence we have  $\max\{m_{1,n},m_{2,n}\}\to +\ii$, $\ainf$.\\

STEP $2$. Let $x_{i,n}$ be any maximum point of $\psi_{i,n}$, $i=1,2$, we prove that
$$
{\rm dist}(x_{i,n},\pa \om)\geq d_0,
$$
for some positive constant $d_0>0$.

Suppose first $N=2$. We argue as in \cite{gnn} p.223.
Let $x_0\in \pa\om$, $\nu_0$ be the outer unit normal at
$x_0$ and $B_r(x_1)\subset \R^N\setminus \ov{\om}$ such that $\ov{B}_r(x_1)\cap \pa\om=\{x_0\}$.
This is always possible since $\om$ is of class $C^3$.
After a translation and a rotation we can assume w.l.o.g. that $x_1=0$, $\nu_0=(-1,0,\cdots,0)$.
Also, after a dilation and a suitable scaling of $\lm_n$, we can assume that $r=1$, whence $x_0=(1,0,\cdots,0)$.\\
At this point if $N=2$ we define (see \cite{gnn} p.223) $y=\frac{x}{|x|^2}$ and $u_{i,n}(y)=\psi_{i,n}(x)$,
$i=1,2$.
The image of $\om$ under this Kelvin
transform, say $\widetilde{\om}$, lies inside $B_1(0)$ and its closure touches the boundary only at $x_0$. Also
$(u_{1,n}(y),u_{2,n}(y))$ satisfies
$$
\graf{-\Delta u_{1,n} =h(y)(\al_{2,n}+{\lm_n}u_{2,n})^{p_2}\quad \mbox{in}\;\;\widetilde{\om}\\ \\
-\Delta u_{2,n} =h(y)(\al_{1,n}+{\lm_n}u_{1,n})^{p_1}\quad \mbox{in}\;\;\widetilde{\om}\\ \\
u_{i,n}=0 \quad \mbox{on}\;\;\pa\widetilde{\om},\quad i=1,2, \\ \\
\al_{i,n}\geq 0,\quad i=1,2,
}
$$
where $h(y)=\frac{1}{|y|^4}$. Let $y=(y_1,y_2)$, clearly $\frac{\pa }{\pa y_1}h(y)<0$ as far as $y_1>0$,
whence we can apply the argument of Theorem 2.1$^{'}$ in \cite{gnn}, just
replacing the classical strong maximum principle and Hopf lemma used there (see Lemma H in \cite{gnn})
with the corresponding results for cooperative and strongly coupled linear elliptic systems, see Theorem 2.2 in
\cite{DP}. Therefore, if $N=2$, as in \cite{gnn} p.223 we deduce that in a neighborhood of $x_0$ depending
only by the geometry of $\om$ there are no critical points of $\psi_{i,n}$, $i=1,2$. Since $\pa \om$ is compact the desired conclusion follows by a covering argument.\\
Suppose now $N\geq 3$ and $\Omega$ convex. As above, using the results for cooperative and strongly coupled linear elliptic systems we can exploit
the argument of Theorem 2.1$^{'}$ in \cite{gnn}. It is well known that the convexity condition ensures then that there are no critical points of $\psi_{i,n}$, $i=1,2$ in a sufficiently small uniform neighborhood of the boundary, see \cite{dFLN} for further details in the scalar case.\\
\\

STEP 3. We obtain a contradiction assuming that $m_{2,n}\leq m_{1,n}\to +\ii$, $\ainf$.
We will never use the fact that $p_1\leq p_2$, whence a contradiction arises in the same
way in the case where, along a subsequence, $m_{1,n}\leq m_{2,n}\to +\ii$, $\ainf$.\\
Let
$x_{1,n}$ be such that $\psi_{1,n}(x_{1,n})=m_{1,n}$. Then by step 2 we have $x_{i,n}\to \ov{x}\in \om$ and after a
translation we can
assume w.l.o.g. that ${x}_{1,n}=0$, $\forall \,n\in\N$.
There are only two possibilities: either,\\
$(j)$ $\sup\limits_{n}\frac{m^{p_2+1}_{2,n}}{m^{p_1+1}_{1,n}}\leq C$, or, passing to a further subsequence if necessary,\\
$(jj)$ $\frac{m^{p_1+1}_{1,n}}{m^{p_2+1}_{2,n}}\to 0$, $\ainf$.\\

We discuss $(j)$ first and define $\dt^2_n=\frac{m_{2,n}}{\lm m^{p_1}_{1,n}}$ and, for $i=1,2$,

$$
v_{i,n}(y)=\frac{1}{m_{i,n}}(\al_{i,n}+\lm \psi_{i,n}(\dt_n y)),\quad y \in \om_{n}=\{y\in \R^N\,:\,\dt_n y\in\om \}
$$
which satisfy
$$
\graf{-\Delta v_{1,n} =\frac{m^{p_2+1}_{2,n}}{m^{p_1+1}_{1,n}}\,v_{2,n}^{p_2}\quad \mbox{in}\;\;\om_{n}\\ \\
-\Delta v_{2,n} =v_{1,n}^{p_1}\quad \mbox{in}\;\;\om_{n}\\ \\
v_{i,n}(y)\leq v_{1,n}(0)= 1 \quad \mbox{in}\;\;\om_n,\quad i=1,2,\\ \\
v_{i,n}(y)\geq \frac{\al_{i,n}}{m_{i,n}}\geq 0 \quad \mbox{in}\;\;\om_n,\quad i=1,2.
}
$$

Since $0\in\om$, then for any $R\geq 1$ we have that for any $n$ large enough it holds $B_{R}(0)\subset \om_n$.
Along a subsequence we can assume that $\frac{m^{p_2+1}_{2,n}}{m^{p_1+1}_{1,n}}\to \mu\in [0,+\ii)$.
Since $\|v_{i,n}\|_{L^{\ii}(\om_n)}\leq C$, then by standard elliptic estimates there exists a subsequence
such that $v_{i,n}$, $i=1,2$ converge in
$C^{2}_{\rm loc}(\R^N)$ to $v_{i}$, $i=1,2$ which are classical solutions of
\beq\label{sys1}
\graf{-\Delta v_{1} =\mu\,v_{2}^{p_2}\quad \mbox{in}\;\;\R^N\\ \\
-\Delta v_{2} =v_{1}^{p_1}\quad \mbox{in}\;\;\R^N\\ \\
0\leq v_{i}(y)\leq v_{1}(0)= 1 \quad \mbox{in}\;\;\R^N,\quad i=1,2,
}
\eeq

At this point observe that if $\mu=0$ then necessarily $v_{1}\equiv 1$ in $\R^N$ and then
$v_2$ would solve,
$$
\graf{
-\Delta v_{2} =1\quad \mbox{in}\;\;\R^N\\ \\
0\leq v_{2}(y)\leq  1 \quad \mbox{in}\;\;\R^N,\quad i=1,2.
}
$$
By the maximum principle we would have $v_2(y)\geq \frac{1}{2N}(R^2-|y|^2)$, for any $R>0$ and in particular
$v_2(0)\geq \frac{1}{2N}R^2$, for any $R>0$, which is impossible. Therefore $\mu\in(0,+\ii)$ which however
is also impossible since it is well known by the result in \cite{Souto} that
\eqref{Souto}
implies that the unique
solution of \rife{sys1} is $v_{i}\equiv 0$, $i=1,2$.\\
Therefore $(jj)$ holds and in this case we choose $\dt^2_n=\frac{m_{1,n}}{\lm m^{p_2}_{2,n}}$ so that
$(v_{1,n},v_{2,n})$ satisfies

$$
\graf{-\Delta v_{1,n} =\,v_{2,n}^{p_2}\quad \mbox{in}\;\;\om_{n}\\ \\
-\Delta v_{2,n} =\frac{m^{p_1+1}_{1,n}}{m^{p_2+1}_{2,n}}v_{1,n}^{p_1}\quad \mbox{in}\;\;\om_{n}\\ \\
v_{i,n}(y)\leq v_{1,n}(0)= 1 \quad \mbox{in}\;\;\om_n,\quad i=1,2,\\ \\
v_{i,n}(y)\geq \frac{\al_{i,n}}{m_{i,n}}\geq 0 \quad \mbox{in}\;\;\om_n,\quad i=1,2.
}
$$
This case is easily seen to lead to the same situation described above for $\mu=0$, whence a contradiction arise
in this case as well.\\
As mentioned above, by symmetry, the
discussion of the case in which along a subsequence $m_{1,n}\leq m_{2,n}\to +\ii$, $\ainf$,
is exactly the same.
\finedim

\bigskip
\bigskip

\section{\bf Spectral and bifurcation analysis}\label{sec2}
In this section we develop the spectral and bifurcation analysis for \un{positive} solutions of $\prl$ with
$\lm\geq 0$ and $(p_1,p_2)$ satisfying \eqref{Souto}. From now on and unless otherwise specified,
$(\bal,\bpl)$ is assumed to be a positive solution of
$\prl$.\\

By the maximum principle $\plv[i]\geq 0$, $i=1,2$ in $\om$ for any solution, whence
for non negative solutions ($\al_{i,\sscp \lm}\geq 0$) we have $\alv[i]+\plv[i]\equiv (\alv[i]+\plv[i])_{+}$.
Therefore from now on and unless otherwise specified we will denote by,
$$
\tlv[i]=\lm p_i,\,\rlav[i](\psi_i)=(\al_i+\lm\psi_i)^{p_i},\,
\rlv[i]=(\alv[i]+\lm\plv[i])^{p_i},\,i=1,2,
$$
$$
\rlavv[i](\psi_i)=(\al_i+\lm\psi_i)^{p_i-1} \mbox{ and } \rlqv[i]=(\alv[i]+\lm\plv[i])^{p_i-1},
$$

where $\alv[i],\plv[i]$, $i=1,2$ denote the components of a non negative solution $(\bal,\bpl)$ of $\prl$ and $q_i$ denotes the conjugate exponent of $p_i$, that is,
$$
\frac1p_i+\frac1q_i=1,\quad i=1,2.
$$

For $(\bal,\bpl)$ a non negative solution of $\prl$ we denote,
$$
<\eta>_{i,\ssl}=\frac{\ino \rlqv[i]  \eta}{\ino \rlqv[i]}\quad \mbox{and }\quad
[\eta]_{i,\ssl}=\eta \,-<\eta>_{i,\ssl},\,i=1,2,
$$
and define,
$$
<\eta,\phi>_{i,\ssl}:={\ino \rlqv[i]  \eta\phi}\quad \mbox{and}\quad
\|\phi\|_{i,\ssl}^2:=<\phi,\phi>_{i,\ssl}={\ino \rlqv[i]  \phi^2}\,,i=1,2,
$$
where $\{\eta,\phi\}\subset L^2(\om)$. For non negative solutions $(\bal,\bpl)$ of $\prl$,
by the strong maximum principle we have that
$\rlv[i]$, $i=1,2$, is strictly positive in ${\om}$, whence $<\cdot,\cdot>_{i,\ssl}$, $i=1,2$
define scalar products on $L^2(\om)$ whose norms are denoted by
$\|\cdot\|_{i,\ssl}$, $i=1,2$.
We will also adopt the useful shorthand notation,
$$
\mlv[i]=\ino \rlqv[i],\,i=1,2.
$$
In the sequel we aim to describe possible branches of solutions of $\prl$ around a positive solution, i.e. with $\alv[i]>0,\,i=1,2$. To this end, it is not difficult to construct an open subset $A_{\sscp \om}$ of the Banach space of triples $(\lm,\ba,\bp)\in \R\times \R^2\times  (C^{2,r_0}_0(\ov{\om}\,))^2$  such that, on $A_\om$, the densities $\rlav[i]=\rlav[i](\psi_i)=(\al_i+\lm\psi_i)^{p_i}$ are well defined and
\begin{equation} \label{bound-a}
\alv[i]+\lm\plv[i]\geq \frac{\alv[i]}{2} \quad \mbox{ in }\ov{\om},\quad i=1,2,
\end{equation}
in a sufficiently small open neighborhood in $A_{\sscp \om}$ of any triple of the form
$(\lm,\bal,\bpl)$ whenever $(\bal,\bpl)$ is a positive ($\al_{i,\sscp \lm}>0$) solution of {\rm $\prl$}.
At this point we introduce the maps,
\beq\label{eF}
\bF: A_{\sscp \om} \to  (C^{r}(\ov{\om}\,))^2,\quad
\bF(\lm,\ba,\bp):=\left(\begin{array}{cl}-\Delta \psi_1 -\rlav[2](\psi_2)\\ \\
-\Delta \psi_2 -\rlav[1](\psi_1)\end{array}\right)
\eeq
and
\beq\label{ePhi}
\Phi:A_{\sscp \om} \to  \R^2\times (C^{r}(\ov{\om}\,))^2,\quad
\Phi(\lm,\ba,\bp):=\left(\begin{array}{cl}\bF(\lm,\ba,\bp)\\-1+\ino \rlav[1](\psi_1)\\
-1+\ino \rlav[2](\psi_2) \end{array}\right),
\eeq
and, for positive solutions and for a fixed $(\lm,\ba,\bp)\in A_{\sscp \om}$,
their differentials with respect to $(\ba,\bp)$, that is the linear operators,
$$
D_{\ba,\bp}\Phi(\lm,\ba,\bp):\R^2 \times (C^{2,r}_0(\ov{\om}\,))^2 \to  \R^2\times (C^{r}(\ov{\om}\,))^2,
$$
which acts as follows on $(\bs,\bphi)=(s_1,s_2,\phi_1,\phi_2)\in \R^2 \times (C^{2,r}_0(\ov{\om}\,))^2$,
\beq\label{dap}
D_{\ba,\bp}\Phi(\lm,\ba,\bp)[\bs,\bphi]=
\left(\begin{array}{cr}
D_\bp \bF(\lm,\ba,\bp)[\bphi]+d_{\ba}\bF(\lm,\ba,\bp) [\bs]\\ \\
\ino \left(D_{\psi_1} \rlav[1](\psi_1)[\phi_1]+d_{\al_1}\rlav[1](\psi_1) [s_1]\right)\\ \\
\ino \left(D_{\psi_2} \rlav[2](\psi_2)[\phi_2]+d_{\al_2}\rlav[2](\psi_2) [s_2]\right)
\end{array}\right),
\eeq
where we have introduced the differential operators,
\beq\label{lin}
D_\bp \bF(\lm,\ba,\bp)[\bphi]=
\left(\begin{array}{cl}
-\Delta \phi_1 -\tlv[2] \rlavv[2](\psi_2) \phi_2\\ \\
-\Delta \phi_2 -\tlv[1] \rlavv[1](\psi_1) \phi_1\end{array}\right),
\quad \bphi=(\phi_1,\phi_2) \in (C^{2,r}_0(\ov{\om}\,))^2,
\eeq

\beq\label{lin1}
D_{\psi_i} \rlav[i][\phi_i]=\tlv[i]V_{i,\lm,\al_{i}}(\psi_i)\phi_i, \quad \phi_i \in C^{2,r}_0(\ov{\om}\,),\,i=1,2,
\eeq

and

\beq\label{lin2}
d_{\ba}\bF(\lm,\ba,\bp) [\bs]=\left(\begin{array}{cl}
-p_2 \rlavv[2](\psi_2) s_2\\ \\
-p_1 \rlavv[1](\psi_1) s_1\end{array}\right)\quad \bs=(s_1,s_2) \in \R^2,
\eeq

\beq\label{lin3}
d_{\al_i}\rlav[i][s_i]= p_i\rlavv[i](\psi_i) s_i, \quad s_i \in \R,\,i=1,2,
\eeq
where we recall $\tlv[i]=\lm p_i$.\\
By the construction of $A_{\sscp \om}$, see in particular \eqref{bound-a}, relying on
known techniques about real analytic functions on Banach spaces {\rm (}\cite{but}{\rm )},
it can be shown that $\Phi(\lm,\ba,\bp)$ is jointly real analytic in an open neighborhood of $A_{\sscp \om}$
around any triple of the form $(\lm,\bal,\bpl)$ whenever $(\bal,\bpl)$ is a positive solution of {\rm $\prl$}.\\

For fixed $\lm> 0$ the pair $(\bal,\bpl)$ is a non negative solution of $\prl$ in the classical sense
as defined in the introduction
if and only if
$\Phi(\lm,\bal,\bpl)=(0,0,0,0)$, and, for positive solutions, we define the linearized operator,
\beq\label{eLl}
{\mathbf L}_{\ssl}[\bphi]=
\left(\begin{array}{cl}
L_{1,\ssl}[\bphi]\\ \\
L_{2,\ssl}[\bphi]
\end{array}\right)=
\left(\begin{array}{cl}
-\Delta \phi_1 -\tlv[2] \rlqv[2] [\phi_2]_{2,\ssl},\\ \\
-\Delta \phi_2 -\tlv[1] \rlqv[1] [\phi_1]_{1,\ssl}\end{array}\right).
\eeq

We say that $\sg=\sg(\bal,\bpl)\in\R$ is an eigenvalue of ${\mathbf L}_{\ssl}$ if the "Hamiltonian" eigenvalue
equation,
\beq\label{lineq0}
\graf{
-\Delta \phi_1-\tlv[2] \rlqv[2] [\phi_2]_{2,\ssl}=\sg p_2\rlqv[2] [\phi_2]_{2,\ssl},\\
-\Delta \phi_2-\tlv[1] \rlqv[1] [\phi_1]_{1,\ssl}=\sg p_1\rlqv[1] [\phi_1]_{1,\ssl},
}
\eeq
admits a non-trivial weak solution $(\phi_1,\phi_2)\in H^1_0(\om)\times H^1_0(\om)$, which is by definition
an eigenfunction of $\sg$.\\
We show that, with this particular definition, the eigenvalues of ${\mathbf L}_{\ssl}$ share the
usual properties of a general self-adjoint elliptic operator.

Let us define the Hilbert space,
\beq\label{y0}
{\mathbf Y}_0:=\left\{ \bvarp=(\varphi_1,\varphi_2) \in \{(L^2(\om))^2,<\cdot,\cdot >_{\ssl}\}
:\,<\varphi_i>_{i,\ssl}=0,i=1,2\right\},
\eeq
where
$$
<\bvarp,{\boldsymbol \eta}>_{\ssl}:=p_1<\varphi_1,\eta_1>_{1,\ssl}+p_2<\varphi_2,\eta_2>_{2,\ssl},
\forall\, (\bvarp,{\boldsymbol \eta})\in ({\mathbf Y}_0)^2
$$
and the linear operator ${\mathbf T}_0:{\mathbf Y}_0\to {\mathbf Y}_0$, which acts on
$\bvarp=(\varphi_1,\varphi_2)\in {\mathbf Y}_0$ as follows
\beq\label{T0}
{\mathbf T}_0(\bvarp):=\left(\begin{array}{cl}
\tlv[2]G[\rlqv[2] \varphi_2]-<\tlv[2]G[\rlqv[2] \varphi_2]>_{1,\ssl}\\ \\
\tlv[1]G[\rlqv[1] \varphi_1]-<\tlv[1]G[\rlqv[1] \varphi_1]>_{2,\ssl}\end{array}\right)\equiv
\left(\begin{array}{cl}
\tlv[2][G[\rlqv[2] \varphi_2]]_{1,\ssl}\\ \\
\tlv[1][G[\rlqv[1] \varphi_1]]_{2,\ssl}\end{array}\right).
\eeq
By standard elliptic theory $G[\rlqv[i] \varphi_i]\subset W^{2,2}(\om)$, $i=1,2$, whence ${\mathbf T}_0$ is compact by the Sobolev embedding.
Also for any $(\bvarp,{\boldsymbol \eta})\in ({\mathbf Y}_0)^2$, we have,
$$
<{\boldsymbol \eta},\mathbf T_0(\bvarp)>_{\ssl}=p_1<\eta_1,[\tlv[2]G[\rlqv[2] \varphi_2]]_{1,\ssl}>_{1,\ssl}+
p_2<\eta_2,[\tlv[1]G[\rlqv[1] \varphi_1]]_{2,\ssl}>_{2,\ssl}=
$$
$$
\lm p_1 p_2\ino \rlqv[1]\eta_1 [G[\rlqv[2] \varphi_2]]_{1,\ssl}+
\lm p_2 p_1\ino \rlqv[2]\eta_2 [G[\rlqv[1] \varphi_1]_{2,\ssl}=
$$
$$
\lm p_1 p_2\ino \rlqv[1]\eta_1 G[\rlqv[2] \varphi_2]+
\lm p_1 p_2\ino \rlqv[2]\eta_2 G[\rlqv[1] \varphi_1]=
$$
$$
\lm p_1 p_2\ino G[\rlqv[1]\eta_1] \rlqv[2] \varphi_2+
\lm p_1 p_2\ino G[\rlqv[2]\eta_2] \rlqv[1] \varphi_1=
$$

$$
p_2\ino [\tlv[1]G[\rlqv[1]\eta_1]]_{2,\ssl} \rlqv[2] \varphi_2+
p_1\ino [\tlv[2]G[\rlqv[2]\eta_2]]_{1,\ssl} \rlqv[1] \varphi_1=
$$

$$
p_1<[\tlv[2]G[\rlqv[2]\eta_2]]_{1,\ssl},\varphi_{1}>_{1,\ssl}+
p_2<[\tlv[1]G[\rlqv[1]\eta_1]]_{2,\ssl},\varphi_{2}>_{2,\ssl}=<{\mathbf T}_0({\boldsymbol \eta}),\bvarp>_{\ssl},
$$
which shows that ${\mathbf T}_0$ is also self-adjoint. Remark also that
$<\bvarp,{\mathbf T}_0(\bvarp)>_{\ssl}>0$ if $\bvarp\neq (0,0)$, as is readily verified observing that
$\psi_i=G[\rlqv[i]\varphi_i]$ satisfies $\ino|\nabla \psi_i|^2=<\varphi_{i},G[\rlqv[i]\varphi_i]>_{i,\ssl}$, $i=1,2$.\\

As a consequence, by the spectral Theorem for
self-adjoint, compact, linear operators on
Hilbert spaces, we have that ${\mathbf Y}_0$ is the Hilbertian direct sum of the eigenfunctions of
${\mathbf T}_0$, which can be represented as follows,
$$
\bvarp_k=(\varphi_{1,k},\,\varphi_{2,k}),\varphi_{i,k}=[\phi_{i,k}]_{i,\ssl},\,i=1,2,\, k\in\N=\{1,2,\cdots\},
$$
$$
{\mathbf Y}_0=\overline{\mbox{Span}\left\{([\phi_{1,k}]_{1,\ssl},[\phi_{2,k}]_{2,\ssl}),\;k\in \N\right\}},
$$
for some $(\phi_{1,k},\phi_{2,k})\in (H^1_0(\om))^2$, $k\in\N=\{1,2,\cdots\}$. In fact, any eigenfunction $\bvarp_k$,
whose eigenvalue is $\mu_k\in \R\setminus\{0\}$, satisfies,
$\mu_k\bvarp_k={\mathbf T}_0(\bvarp_k)$ and, by defining,
$$
\frac{\lm}{\lm+\sg_k}=\mu_k\in (0,+\ii),\quad 0>\mu_1\geq \mu_2\geq....\geq\mu_k\to0,
$$
$$
\phi_{1,k}:=(\tlv[2]+p_2\sg_k) G[\rlqv[2]\varphi_{2,k}],\,
\phi_{2,k}:=(\tlv[1]+p_1\sg_k) G[\rlqv[1]\varphi_{1,k}],
$$
it is easy to see that $\bvarp_k$ is an eigenfunction of ${\mathbf T}_0$ with eigenvalue $\mu_k$ if and only if
$\bphi_k=(\phi_{1,k},\phi_{2,k})\in (H^1_0(\om))^2$ weakly solves,
\beq\label{lineqm}
\graf{
-\Delta \phi_{1,k}=(\tlv[2]+\sg_k p_2)\rlqv[2] [\phi_{2,k}]_{2,\ssl},\\
-\Delta \phi_{2,k}=(\tlv[1]+\sg_k p_1)\rlqv[1] [\phi_{1,k}]_{1,\ssl}.}
\eeq

In particular the first eigenvalue $\sg_1=\sg_1(\bal,\bpl)$ is uniquely defined by the spectral radius of ${\mathbf T}_0$, $r({\mathbf T}_0)\equiv \mu_1=\frac{\lm}{\lm+\sg_1}$ where
$$
\mu_1=r({\mathbf T}_0)= \sup\limits_{\bvarp\in{\mathbf Y}_0\setminus \{\mathbf 0\}}
\frac{<\bvarp,T_0(\bvarp)>_{\ssl}}{<\bvarp,\bvarp>_{\ssl}}.
$$
Since $\sg_1+\lm=\frac{\lm}{\mu_1}$ and since $\mu_1$ is positive, then
\beq\label{lineq2}
\lm+\sg_1>0.
\eeq

By the Fredholm alternative,
if $0\notin \{\sg_j\}_{j\in \N}$, then ${\mathbf I}-{\mathbf T}_0$ is an isomorphism of ${\mathbf Y}_0$ onto itself. Clearly, we can construct an orthonormal base of
eigenfunctions $\{\bphi_k\}_{k\in\N}$ with respect to the scalar product $<\cdot,\cdot>_{\ssl}$. However we need a refined property which is the following,
\ble\label{ortbase}
There exists a complete orthonormal base  $\{{\boldsymbol \varphi}_k\}_{k\in\N}$ of ${\boldsymbol Y}_0$,
with the property that ${\varphi}_k=([\phi_{1,k}]_{1,\ssl},[\phi_{2,k}]_{2,\ssl})$ satisfies,
\beq\label{orto}
<[\phi_{1,k}]_{1,\ssl}, [\phi_{1,j}]_{1,\ssl}>_{1,\ssl}=0=
<[\phi_{2,k}]_{2,\ssl}, [\phi_{2,j}]_{2,\ssl}>_{2,\ssl},\quad \forall k\neq j.
\eeq
In particular $\{[\phi_{i,k}]_{}\}_{k\in \N}$, $i=1,2$ is a complete orthonormal base of
$$
Y_{i,0}:=\left\{ \varphi_i \in \{L^2(\om),<\cdot,\cdot >_{i,\ssl}\} :\,<\varphi_i>_{i,\ssl}=0\right\},
i=1,2.
$$
\ele
\proof
If $\sg_k,\sg_j$ denote the eigenvalues of ${\boldsymbol \phi}_k$, ${\boldsymbol \phi}_j$ respectively,
and if $\sg_k\neq\sg_j$, then for fixed $k$ we can multiply the first equation in \rife{lineqm} by
$\phi_{2,j}$, the second by $\phi_{1,j}$ and integrate by parts to conclude that
$$
\graf{
(\tlv[1]+\sg_j p_1)\ino \rlqv[1] [\phi_{1,j}]_{1,\ssl}[\phi_{1,k}]_{1,\ssl}=
(\tlv[2]+\sg_k p_2)\ino \rlqv[2] [\phi_{2,k}]_{2,\ssl}[\phi_{2,j}]_{2,\ssl},\\
(\tlv[2]+\sg_j p_2)\ino \rlqv[2] [\phi_{2,j}]_{2,\ssl}[\phi_{2,k}]_{2,\ssl}=
(\tlv[1]+\sg_k p_1)\ino \rlqv[1] [\phi_{1,k}]_{1,\ssl}[\phi_{1,j}]_{1,\ssl}.}
$$
Putting
$$
x=p_1\ino \rlqv[1] [\phi_{1,j}]_{1,\ssl}[\phi_{1,k}]_{1,\ssl},\; y=p_2 \ino \rlqv[2] [\phi_{2,k}]_{2,\ssl}[\phi_{2,j}]_{2,\ssl},
$$
this is equivalent to the system
$$
\graf{(\lm+\sg_j)x-(\lm+\sg_k)y=0\\ -(\lm+\sg_k)x+(\lm+\sg_j)y=0},
$$
whose determinant is not zero as far as $\sg_k\neq \sg_j$. Therefore in particular
$$
<[\phi_{i,j}]_{i,\ssl}[\phi_{i,k}]_{i,\ssl}>_{i,\ssl}=\ino \rlqv[i] [\phi_{i,j}]_{1,\ssl}[\phi_{i,k}]_{i,\ssl}=0,\,i=1,2,
$$
whenever $\sg_j\neq\sg_k$. If $\sg_k=\sg_j$, since the eigenspace is of finite dimension, a
standard componentwise orthonormalization
argument shows that in fact the basis can be chosen to satisfy \rife{orto}.\\
At this point, since $\{\bvarp_k\}_{k\in\N}$ is a complete orthonormal base then
also $\{[\phi_{i,k}]_{i,\ssl}\}_{k\in \N}$ must be a
complete orthonormal base of $Y_{i,0}$ for any $i=1,2$.
\finedim

\bigskip
\bigskip

Concerning $D_{\ba,\bp}\Phi(\lm,\ba,\bp)$ we have,

\bpr\label{pr2.2} For any positive solution $(\bal,\bpl)$ of {\rm$\prl$} with $\lm\geq 0$, the kernel of
$D_{\ba,\bp}\Phi(\lm,\bal,\bpl)$ is empty if and only if the system,
{\rm
\beq\label{contr}
\graf{
-\Delta \phi_1-\tlv[2] \rlqv[2] [\phi_2]_{2,\ssl}=0,\\
-\Delta \phi_2-\tlv[1] \rlqv[1] [\phi_1]_{1,\ssl}=0,
}\quad (\phi_1,\phi_2)\in (C^{2,r}_0(\ov{\om}\,))^2
\eeq
}
admits only the trivial solution, or equivalently, if and only if $0$ is not an
eigenvalue of ${\mathbf L}_{\ssl}$.
\epr
\proof
If $(\phi_1,\phi_2)\in (H^{1}_0(\om))^2$ solves \rife{contr} and since $\om$ is of class $C^3$, then by standard elliptic regularity theory
and a bootstrap argument we have
$(\phi_1,\phi_2)\in (C^{2,r}_0(\ov{\om}\,))^2$. Therefore, in particular $0$ is not an eigenvalue of ${\mathbf L}_{\ssl}$ if and only if \rife{contr}
admits only the trivial solution.\\
Suppose first that there exists a non-vanishing pair $(\bs,\bphi)\in \R^2 \times (C^{2,r}_0(\ov{\om}\,))^2 $ such that
$$D_{\ba,\bp}\Phi(\lm,\bal,\bpl)[\bs,\bphi]=({\bf 0},{\bf 0}).$$
Then the second pair of equations in \rife{dap}, that is
$\left.\ino \left(D_\psi \rlv[i][\phi_i]+d_{\al_i}\rlv[i] [s_i]\right)\right|_{(\ba,\bp)=(\bal,\bpl)}=0$, take the form,
$$
p_i \ino \left(\lm\rlqv[i] \phi_i+\rlqv[i] s_i\right)=0,\,i=1,2
$$
which is equivalent to
$$
s_i=s_{i,\ssl}=- \lm<\phi_i>_{i,\ssl}.
$$
Substituting this relations into the first pair of equations,
$$
{\bf L}_{\ssl}[\bphi]=D_\bp F(\lm,\bal,\bpl)[\bphi]+d_{\ba}F(\lm,\bal,\bpl) [\bs_{\ssl}]=0,
$$
we conclude that $\bphi=(\phi_1,\phi_2)$ is a non-trivial, classical
solution of \rife{contr}.\\ This shows one part of the claim, while on the other side, if a non-trivial, classical
solution of \rife{contr} exists, then by arguing the other way around, obviously we can find some $(\bs,\bphi)\neq ({\bf 0},{\bf 0})$
such that $D_{\ba,\bp}\Phi(\lm,\bal,\bpl)[\bs,\bphi]=({\bf 0},{\bf 0})$, as claimed.
\finedim

\bigskip
\bigskip

We state now the result needed to describe branches of solutions
of $\prl$ at regular points.
\ble\label{lem1.1} Let $(\ba_{\sscp \lm_0},\bp_{\sscp \lm_0})$ be a positive solution of {\rm $\prl$} with $\lm=\lm_0\geq 0$.\\
If $0$ is not an eigenvalue of ${\bf L}_{\sscp \lm_0}$, then:\\
$(i)$ $D_{\ba,\bp}\Phi(\lm_0, \ba_{\sscp \lm_0}, \bp_{\sscp \lm_0})$ is an isomorphism;\\
$(ii)$ There exists an open neighborhood $\mathcal{U}\subset A_{\sscp \om}$ of $(\lm_0,\ba_{\sscp \lm_0},\bp_{\sscp \lm_0})$ such
that the set of solutions of
{\rm $\prl$} in $\mathcal{U}$ is a real analytic curve of positive solutions $J\ni\lm\mapsto (\bal,\bpl)\in B$, for
suitable neighborhoods $J$ of $\lm_0$ and $B$ of $(\ba_{\sscp \lm_0},\bp_{\sscp \lm_0})$ in
$(0,+\ii)^2\times (C^{2,r}_{0,+}(\ov{\om}\,))^2$.\\
$(iii)$ In particular if $(\ba_{\sscp \lm_0},\bp_{\sscp \lm_0})=(\ba_{0},\bp_{0})$, then
$(\bal,\bpl)=(\ba_0,\bp_{0})+\mbox{\rm O}(\lm)$ as $\lm\to 0$.
\ele

\proof
By the construction of $A_{\sscp \om}$, the map
$F$ as defined in \rife{eF} is jointly analytic in a suitable neighborhood of $(\lm,\bal,\bpl)$. As a consequence, whenever $(i)$ holds,
then $(ii)$ is an immediate consequence of the real analytic implicit function theorem,
see for example Theorem 4.5.4 in \cite{but}. In particular $(iii)$ is a straightforward consequence of $(ii)$. Therefore, we are
just left with the proof of $(i)$.\\
Concerning $(i)$ we observe that, although the differential of the constrained equations
(which are the last two differentials in \rife{dap}), do not define a Fredholm operator (since obviously the dimension of their kernel is not finite dimensional), however a simple inspection shows that
in fact $D_{\ba,\bp}\Phi(\lm_0, \ba_{\sscp \lm_0}, \bp_{\sscp \lm_0})$ is a Fredholm operator, see for example Lemma 2.4 in \cite{BJ1}.
As a consequence $(i)$ follows from Proposition \ref{pr2.2} and the Fredholm alternative.\finedim

\bigskip
\bigskip

We conclude this section with some spectral estimates about $\sg_1(\bal,\bpl)$ and $\lm^*(\om,{\mathbf p})$.
\bpr\label{prsigma} Let $(p_1,p_2)$ satisfy \eqref{Souto}, assume without loss of generality $p_1\leq p_2$ and suppose that $(\bal,\bpl)$ is
a positive solution of {\rm $\prl$} with $\lm\leq \frac{1}{p_2}\Lambda(\om,2p_2)$. Then $\sg_1(\bal,\bpl)> 0$.
\epr
\proof
For $k=1$ and to ease the notations, in this proof we set $(\phi_1,\phi_2)=(\phi_{1,1},\phi_{2,1})$, where $(\phi_{1,1},\phi_{2,1})$ is any eigenvector of $\sg_1$. Clearly we can just consider $\lambda>0$. We multiply the first equation in \rife{lineqm} by $\phi_{1}$, the second by $\phi_{2}$ and deduce that
$$
\ino |\nabla \phi_1|^2=(\lm+\sg_1)p_2\ino \rlqv[2][\phi_2]_{2,\ssl}[\phi_1]_{1,\ssl},
$$
and,
$$
\ino |\nabla \phi_2|^2=(\lm+\sg_1)p_1\ino \rlqv[1][\phi_1]_{1,\ssl}[\phi_2]_{2,\ssl},
$$
which in view of \rife{lineq2} shows that
\beq\label{den1}
\ino \rlqv[2][\phi_2]_{2,\ssl}[\phi_1]_{1,\ssl}>0,
\ino \rlqv[1][\phi_1]_{1,\ssl}[\phi_2]_{2,\ssl}>0.
\eeq

In particular  we readily deduce that
$$
\sg_1=\frac{N(\phi_1,\phi_2)}
{\ino (p_1\rlqv[1]+p_2\rlqv[2])[\phi_1]_{1,\ssl}[\phi_2]_{2,\ssl}},
$$
where
$$
N(\phi_1,\phi_2)=\ino |\nabla \phi_{1}|^2+\ino |\nabla \phi_{2}|^2-
\lm p_2\ino \rlqv[2][\phi_2]_{2,\ssl}[\phi_1]_{1,\ssl}-\lm p_1\ino \rlqv[1][\phi_2]_{2,\ssl}[\phi_1]_{1,\ssl}.
$$

By using $ab\leq \frac12(a^2+b^2)$ we see that,
\beq\label{nphi}
N(\phi_1,\phi_2)\geq \frac12\sum\limits_{i=1}^2
\left(\ino |\nabla \phi_{i}|^2-\lm p_1\ino \rlqv[1][\phi_i]^2_{i,\ssl}\right)+
\eeq
$$
\frac12\sum\limits_{i=1}^2\left( \ino |\nabla \phi_{i}|^2-\lm p_2\ino \rlqv[2][\phi_i]^2_{i,\ssl}\right)\geq
$$

$$
\frac12\sum\limits_{i=1}^2
\left(\ino |\nabla \phi_{i}|^2-\lm p_1\ino \rlqv[1]\phi_i\right)+
\frac12\sum\limits_{i=1}^2\left( \ino |\nabla \phi_{i}|^2-\lm p_2\ino \rlqv[2]\phi_i\right),
$$
where the last inequality is an equality if and only if $<\phi_1>_{1,\ssl}=0=<\phi_2>_{2,\ssl}$.\\
Next observe that for any $\phi \in H^1_0(\om)\setminus\{0\}$ and for any $\{i,j\}\in \{1,2\}$, by the Holder
inequality and \rife{sob} we have,
\beq\label{firste1}
\ino |\nabla \phi|^2-\lm p_i\ino \rlqv[i]\phi^2\geq
\ino |\nabla \phi|^2-\lm p_i\left(\ino \rlv[i]\right)^{\frac{1}{q_i}}\left(\ino \phi^{2p_i}\right)^{\frac{1}{p_i}}=
\eeq
\beq\label{firste2}
\ino |\nabla \phi|^2-\lm p_i\left(\ino \phi^{2p_i}\right)^{\frac{1}{p_i}}\geq
\left(\ino \phi^{2p_i}\right)^{\frac{1}{p_i}}\left(\Lambda(\om,2p_i)-\lm p_i\right),\,i=1,2.
\eeq

Since $|\om|=1$ and $p_1\leq p_2$, it is well known (see Theorem 3 in \cite{CRa})
that $\Lambda(\om,2p_2)\leq \Lambda(\om,2p_1)$ and that the inequality is strict as far as $p_1<p_2$.
Therefore we have proved that if $\lm\leq \frac1{p_2}\Lambda(\om,2p_2)$ then $\sg_1\geq 0$.\\
At this point we argue by contradiction and assume that $\sg_1=0$ for some
$\lm\leq \frac{1}{p_2} \Lambda(\om,2p_2)$. Following the equality sign in all the inequalities used so far we see that
 if $\sg_1=0$ then we would necessarily have $\lm=\frac{1}{p_2} \Lambda(\om,2p_2)$, $p_1=p_2$, $<\phi_1>_{1,\ssl}=0=<\phi_2>_{2,\ssl}$,
and in particular $\phi_1=\phi_2$ a.e. in $\om$. By \rife{contr} this implies also
$\rlqv[1]=\rlqv[2]$ a.e. in $\om$. As a consequence we would also have,
$$
0=p_1\sg_1=\frac{\ino |\nabla \phi_1|^2-\lm p_1\ino \rlqv[1]\phi_1^2 }
{\ino \rlqv[1]\phi_1^2}>\nu_1:=
\inf \limits_{\phi\in H^1_0(\om)}\frac{\ino |\nabla \phi|^2-\lm p_1\ino \rlqv[1]\phi^2}
{\ino\rlqv[1]\phi^2},
$$
where the strict inequality is due to the fact that
any function which attains the inf is a first eigenfunction and consequently does not change sign.
On the other side if $\phi \in C^{1}_0(\ov{\om}\,)$, $\phi\equiv \!\!\!\!\!\!/ \;0$, then,
once more by the Holder inequality and \rife{sob}, we have,
$$
\dfrac{\ino |\nabla \phi|^2}{\ino \rlqv[2] \phi ^2}\geq
\dfrac{1}{ \left(\ino \rlv[2]\right)^{\frac{1}{q_2}}}
\dfrac{\ino |\nabla \phi|^2}{\left(\ino \phi^{2p_2}\right)^{\frac{1}{p_2}}}=
\dfrac{\ino |\nabla \phi|^2}{\left(\ino \phi^{2p_1}\right)^{\frac{1}{p_1}}}\geq {\Lambda(\om,2p_1)}
$$
which immediately implies that,
$$
\nu_{1}\geq {\Lambda(\om,2p_1)}-\lm p_1\geq 0,\;\forall\,\lm p_1\leq\Lambda(\om,2p_1),
$$
which is a contradiction to $\nu_{1}<0$.
\finedim

\bigskip

\bpr\label{prsigma1} Let $(p_1,p_2)$ satisfy \eqref{Souto}, assume without loss of generality $p_1\leq p_2$ and suppose that
$\al_{i,*}=0$, $i=1,2$. Then
$\lm^*(\om,{\mathbf p})\geq \frac{1}{p_2}\Lambda(\om,2p_2).$
\epr
\proof
We can assume w.l.o.g. $\lm^*(\om,{\mathbf p})\in(0,+\infty)$. By definition there exists a sequence $\lm_n\to (\lm^*(\om,{\mathbf p}))^{-}$, such
that $\al_{i,\sscp \lm_n}\to 0^+$. By Lemma \ref{lemE1} and passing to a further subsequence if necessary we can assume that
$u_{i,n}=\lm_n\psi_{i,\sscp \lm_n}$, $i=1,2$, converge smoothly to $u_i$, $i=1,2$, which are classical solutions
of
$$
\graf{-\Delta u_{1} =\lm^* u_{2}^{p_2}\quad \mbox{in}\;\;{\om}\\ \\
-\Delta u_{2} =\lm^*u_{1}^{p_1}\quad \mbox{in}\;\;{\om}\\ \\
u_{i}\geq 0 \quad \mbox{in}\;\;\om,\quad i=1,2\\ \\
u_{i}=0 \quad \mbox{on}\;\;\pa \om,\quad i=1,2.
}
$$
By the monotonicity of the energy $\el$ we have $E_{\sscp \lm^*}\geq E_0(\om)$, whence both $u_1$ and $u_2$
cannot vanish identically. Let us set $V_i=u_i^{p_i-1}$, $\phi_i=u_i$, $i=1,2$, then we have,
$$
\graf{-\Delta \phi_{1} =\lm^* V_2\phi_{2}\quad \mbox{in}\;\;{\om}\\ \\
-\Delta \phi_{2} =\lm^*V_1\phi_1 \quad \mbox{in}\;\;{\om}\\ \\
\phi_{i}\geq 0 \quad \mbox{on}\;\;\om,\quad i=1,2\\ \\
\phi_{i}=0 \quad \mbox{on}\;\;\pa \om,\quad i=1,2.
}
$$
By the strong maximum principle for cooperative and strongly coupled linear elliptic systems
(see Theorem 2.2 in \cite{DP}) we have $\phi_i>0$ in $\om$, $i=1,2$ and we deduce that,
$$
\graf{-\Delta \phi_{1} \leq \lm^* p_2V_2\phi_{2}\quad \mbox{in}\;\;{\om}\\ \\
-\Delta \phi_{2} \leq \lm^*p_1 V_1\phi_1 \quad \mbox{in}\;\;{\om}\\ \\
\phi_{i}=0 \quad \mbox{on}\;\;\pa \om,\quad i=1,2,
}
$$
where the equality holds if and only if $p_1=1=p_2$. Multiplying the first equation by
$\phi_1$, the second by $\phi_2$ and integrating by parts we deduce that

\beq\label{2901.1}
0\geq Q(\phi_1,\phi_2):=\ino |\nabla \phi_{1}|^2+\ino |\nabla \phi_{2}|^2-
\lm^* p_2\ino V_2\phi_2\phi_1-\lm^* p_1\ino V_1\phi_2\phi_1.
\eeq

By using $ab\leq \frac12(a^2+b^2)$ we see that,
$$
Q(\phi_1,\phi_2)\geq \frac12\sum\limits_{i=1}^2
\left(\ino |\nabla \phi_{i}|^2-\lm^* p_1\ino V_1\phi_i^2\right)+
\frac12\sum\limits_{i=1}^2\left( \ino |\nabla \phi_{i}|^2-\lm^* p_2\ino V_2\phi_i^2\right),
$$
and then by using \rife{firste1}, \rife{firste2}, $|\om|=1$, $p_1\leq p_2$ and
$\Lambda(\om,2p_2)\leq \Lambda(\om,2p_1)$ (this is well known, see for example \cite{CRa}) as above we deduce that
that
$$
Q(\phi_1,\phi_2)\geq
\left(\frac12\sum\limits_{i=1}^2\ino \ino V_1\phi_i^2+\frac12\sum\limits_{i=1}^2\ino \ino V_2\phi_i^2\right)
({\Lambda(\om,2p_2)}-\lm^* p_2),
$$
and we see from \rife{2901.1} that $\lm^* p_2 \geq {\Lambda(\om,2p_2)}$, as claimed. \finedim

\bigskip
\bigskip

\section{{\bf Existence, uniqueness and monotonicity}}\label{appD}

Let $\mathcal{G}(\om)$ denote the set of solutions of $\prl$ for $\lm\in [0,\lm^*(\om,{\mathbf p}))$.

\ble\label{unique} Let $(p_1,p_2)$ satisfy \eqref{Souto}. For any $\lm\in [0,\lm^*(\om,{\mathbf p}))$ there exists
one and only one solution $(\bal,\bpl)$ of {\rm $\prl$}. In particular $\mathcal{G}(\om)$ is a real analytic simple curve of positive solutions
$[0,\lm^*(\om,{\mathbf p}))\ni \lm\mapsto (\bal,\bpl)$.
\ele
\proof
By Lemma \ref{lmsmall} and Proposition \ref{prsigma} we have that $\lm^*(\om,{\mathbf p})>0$ and
there exists a unique solution of $\prl$ for $\lm\leq \lm_0$ for some $\lm_0$ small enough.
In particular these unique solutions are positive and we can assume, possibly taking a smaller $\lm_0$,
that $\lm^*(\om,{\mathbf p})>\lm_0$. However by definition we have $\sg_1(\all,\pl)>0$ for any
$\lm\in [0,\lm^*(\om,{\mathbf p}))$ and then by Lemma \ref{lemE1} and
Lemma \ref{lem1.1} any solution whose $\lm$ is less than $\lm_0$ generates a real analytic curve
of positive solutions which can be continued to a real analytic curve of positive solutions defined in
$[(-\dt,\lm^*(\om,{\mathbf p}))$ for some small $\dt>0$.\\
If at any point in $(\lm_0,\lm^*(\om,{\mathbf p}))$ there exist two solutions, then by definition of
$\lm^*(\om,{\mathbf p})$ they would be both positive and
each one would generate in the same way another curve
of positive solutions in $(-\dt,\lm^*(\om,{\mathbf p}))$, for some small $\dt>0$.
Obviously for both curves at $\lm=0$ we have $(\all,\pl)=(\al_{\sscp 0},\psi_{\sscp 0})$ which is the
unique solution of $\prl$ for $\lm=0$. This is obviously impossible since then $(\al_{\sscp 0},\psi_{\sscp 0})$ would be a bifurcation
point, in contradiction with Lemma \ref{lem1.1}.
\finedim

\bigskip

Problem $\prl$ is the Euler-Lagrange equation of the constrained
minimization principle {\bf (VP)} below for the densities {\rm $(\rh_1,\rh_2)$}.  As far as $1\leq p_1\leq p_2<p_{N}$,
existence of solutions could be proved by an adaptation of an argument in \cite{BeBr}, worked out there for a more general "scalar" problem,
based on the theory of conjugate convex function. We adopt here a different argument based on the weak Young inequality
(\cite{Lieb}),
which yields existence for $(p_1,p_2)$ satisfying \eqref{Souto}, that is,
\beq\label{Souto-4}
\frac{1}{p_1+1}+\frac{1}{p_2+1}>\frac{N-2}{N-1},\quad p_i\in (0,+\ii),\,i=1,2,
\eeq
whose relevance for elliptic systems of Lane-Emden type was first noticed in \cite{Souto}.

\bte\label{thm5.2} Let $(p_1,p_2)$ satisfy \eqref{Souto}. For any $\lm>0$ there exists a solution $(\bal,\bpl)$ of {\rm $\prl$}.
\ete
\proof
We discuss only the case $N\geq 3$, the case $N=2$ is easier.
Here $\|\rh\|_{p}$ denotes the standard $L^p(\om)$ norm. We will denote with $C$ various constants depending only by
 $(p_1,p_2)$, $\om$ and $N$.\\

Let us define
$r_i=1+\frac{1}{p_i}$, $i=1,2$ and
$$
\mathcal{P}_{\sscp \om,i}:=\left\{\rh\in L^{r_i}(\om)\,|\,\rh\geq 0\;\mbox{a.e. in}\;\om,\;\ino \rh=1 \right\},\,i=1,2.
$$
It is readily seen that \eqref{Souto} is equivalent to
\beq\label{Souto-r}
\frac{1}{r_1}+\frac{1}{r_2}<\frac{N}{N-1}, \quad r_i\in (1,+\ii),\,i=1,2,
\eeq
and for any $(r_1,r_2)$ satisfying \eqref{Souto-r}, for any $\lm>0$ and
$$
(\rh_1,\rh_2)\in \mathcal{P}_{\sscp \om}\equiv \mathcal{P}_{\sscp \om,1}\times \mathcal{P}_{\sscp \om,2},
$$
we define the free energy,
\beq\label{jeil}
J_{\ssl}(\rh_1,\rh_2)=
{\scriptstyle \frac{1}{r_1}}\ino (\rh_1)^{r_1}+
{\scriptstyle \frac{1}{r_2}}\ino (\rh_2)^{r_2}
-\lm\ino \rh_1 G[\rh_2].
\eeq

We split the proof in several steps.\\
STEP 1. We prove that if $(r_1,r_2)$ satisfies \eqref{Souto-r}, then $J_{\ssl}$ is coercive for any $\lm>0$, that is, if $\lm>0$ and $(\rh_{1,n},\rh_{2,n})\in\mathcal{P}_{\sscp \om,1}\times \mathcal{P}_{\sscp \om,2}$ and $\max\{\|\rh_{1,n}\|_{r_1},\|\rh_{2,n}\|_{r_2}\}\to +\ii$, then $J_{\ssl}(\rh_{1,n},\rh_{2,n})\to +\ii$, $\ainf$.\\
Let us recall the weak Young inequality (\cite{Lieb}),
\beq\label{HLS}
\ino \rh_1 G[\rh_2]\leq C\|\rh_1\|_{s_1}\|\rh_2\|_{s_2},\quad \frac{1}{s_1}+\frac{1}{s_2}=\frac{N+2}{N},\quad s_i\in (0,+\ii).
\eeq
Assume for the time being that for any $(r_1,r_2)$ which satisfies \eqref{Souto-r}, there exists $(s_1,s_2)$ which satisfies \eqref{HLS} and
\beq\label{s12.0}
s_i\in (1,r_i),\;i=1,2.
\eeq
Then by the \eqref{HLS} and standard interpolation inequalities (recall $\ino \rh_i=1$, $i=1.2$) we would have that,
\beq\label{HLS.1}
\ino \rh_1 G[\rh_2]\leq C\|\rh_1\|_{s_1}\|\rh_2\|_{s_2}\leq C\|\rh_1\|^{\gamma_1 r_1}_{r_1}\|\rh_2\|^{\gamma_2 r_2}_{r_2},
\eeq
where
$$
\gamma_i=(1-\frac{1}{s_i})\frac{1}{r_i-1},\quad i=1,2.
$$
Based on \eqref{HLS.1}, elementary arguments show that $J_{\ssl}$ is coercive for any $\lm>0$ as far as $\gamma_1+\gamma_2<1$. Therefore
we are left with showing that for any $(r_1,r_2)$ which satisfies \eqref{Souto-r}, there exists $(s_1,s_2)$ which satisfies \eqref{HLS}, \eqref{s12.0} and
in particular,
\beq\label{s12.1}
(1-\frac{1}{s_1})\frac{1}{r_1-1}+(1-\frac{1}{s_2})\frac{1}{r_2-1}<1.
\eeq
Obviously there is no loss of generality in assuming
$$
r_2\geq r_1.
$$
Observe that, as far as,
$$
(r_1,r_2)\in (1,\frac{N}{2}]\times [\frac{N}{2},+\ii)\bigcup  \{ r_2\geq r_1\,:\,[\frac{N}{2},+\ii)\times [\frac{N}{2},+\ii)\},
$$
then putting $\frac{1}{s_2}=\frac{2}{N}+\eps$ for some small enough $\eps>0$, from \eqref{HLS} we would have $\frac{1}{s_1}=1-\eps$ and then
$$
1<s_2<\frac{N}{2}\leq r_2, \quad 1<s_1=\frac{1}{1-\eps}<r_1.
$$
Thus \eqref{s12.0} is satisfied. On the other side, since here $r_2\geq \frac{N}{2}$, we also have that,
$$
(1-\frac{1}{s_1})\frac{1}{r_1-1}+(1-\frac{1}{s_2})\frac{1}{r_2-1}=\frac{\eps}{r_1-1}+(1-\frac{2}{N}-\eps)\frac{1}{r_2-1}=
$$
$$
\eps(\frac{1}{r_1-1}-\frac{1}{r_2-1})+\frac{N-2}{N}\frac{1}{r_2-1}\leq
 \eps(\frac{1}{r_1-1}-\frac{1}{r_2-1})+\frac{2}{N}<1,
$$
for any $\eps>0$ small enough, showing that \eqref{s12.1} is satisfied as well.\\
Therefore we are left with showing that $(s_1,s_2)$ which satisfies \eqref{HLS}, \eqref{s12.0} and \eqref{s12.1} exists in the region
\beq\label{reg-d}
\frac{1}{r_1}+\frac{1}{r_2}<\frac{N}{N-1},\quad r_1\in (r^*_{_N},\frac{N}{2}),\quad r_1\leq r_2\leq \frac{N}{2},
\eeq
where $r^*_{_N}=\frac{N(N-1)}{N^2-N+2}$ is the intersection of $\frac{1}{r_1}+\frac{1}{r_2}<\frac{N}{N-1}$ with the line $r_2=\frac{N}{2}$.
Remark that the lowest possible value of $r_2$ in this region is the one at the corner point where $r_2=r_1=2\frac{N-1}{N}$. Put
$$
\frac{1}{s_1}=\frac{1}{r_1}+\eps,
$$
which in view of \eqref{HLS} implies that
$$
1-\frac{1}{s_2}=\frac{1}{r_1}-\frac{2}{N}+\eps.
$$
Then, in view of \eqref{Souto-r}, and since $r_2-1\geq r_1-1\geq 2\frac{N-1}{N}-1=\frac{N-2}{N} $, concerning \eqref{s12.1} we have,
\begin{align*}
&(1-\frac{1}{s_1})\frac{1}{r_1-1}+(1-\frac{1}{s_2})\frac{1}{r_2-1}=\\
&\eps(\frac{1}{r_2-1}-\frac{1}{r_1-1})+\frac{1}{r_1}+(\frac{1}{r_1}-\frac{2}{N})\frac{1}{r_2-1}\leq\\
&\frac{1}{r_1}+(\frac{1}{r_1}+\frac{1}{r_2}-1-\frac{2}{N})\frac{1}{r_2-1}+(1-\frac{1}{r_2})\frac{1}{r_2-1}=\\
&\frac{1}{r_1}+\frac{1}{r_2}+(\frac{1}{r_1}+\frac{1}{r_2}-1-\frac{2}{N})\frac{1}{r_2-1}<\\
&\frac{1}{r_1}+\frac{1}{r_2}+(\frac{1}{N-1}-\frac{2}{N})\frac{1}{r_2-1}\leq\\
&\frac{1}{r_1}+\frac{1}{r_2}+\frac{2-N}{N(N-1)}\frac{N}{N-2}=\\
&\frac{1}{r_1}+\frac{1}{r_2}-\frac{1}{N-1}<1.
\end{align*}
Therefore we are left with showing that if $(r_1,r_2)$ satisfies \eqref{reg-d}, then we can find $(s_1,s_2)$ which satisfies
\eqref{s12.0} as well as $\frac{1}{s_1}=\frac{1}{r_1}+\eps$, $1-\frac{1}{s_2}=\frac{1}{r_1}-\frac{2}{N}+\eps$ for some $\eps>0$.
We split the discussion in three regions which could be possibly empty for some $N\geq 3$.\\
We recall that the symmetric point of the critical hyperbola in \eqref{HLS} is just $s_1=\frac{2N}{N+2}=s_2$ and that
\beq\label{hyp2}
\frac{1}{r_1}+\frac{1}{r_2}<\frac{N}{N-1}<\frac{N+2}{N}.
\eeq
We start with the domain,
$$
\om_1=\{\,r_1\in [2\frac{N-1}{N},\frac{N}{2}],\,r_1\leq r_2\leq\frac{N}{2}\}.
$$
If $(r_1,r_2)\in\om_1$, since $r_1\in [2\frac{N-1}{N},\frac{N}{2}]$, then for some $\sg_\eps\to 0$, as $\eps\to 0^+$,
we have that if $\frac{1}{s_1}=\frac{1}{r_1}+\eps$ then $r_1>s_1>2\frac{N-1}{N}-\sg_\eps>\frac{2N}{N+2}$, whence in particular
$$
s_2<\frac{2N}{N+2}<2\frac{N-1}{N}\leq r_2,
$$
where we use \eqref{hyp2}. Therefore \eqref{s12.0} is satisfied and then we are done with $\om_1$. Next we consider the case,
$$
\om_2=\{\frac{1}{r_1}+\frac{1}{r_2}<\frac{N}{N-1},\,r_1\in [\frac{2N}{N+2},2\frac{N-1}{N}),\,r_2\leq\frac{N}{2}\}.
$$
If $(r_1,r_2)\in\om_2$, since $r_1\in [\frac{2N}{N+2},2\frac{N-1}{N})$, then for some $\sg_\eps\to 0$, $h_\eps\to 0$, as $\eps\to 0^+$,
we have that if $\frac{1}{s_1}=\frac{1}{r_1}+\eps$ then $r_1>s_1\geq r_1-\sg_\eps 2\frac{2N}{N+2}-\sg_\eps$, whence in particular
$$
s_2<\frac{2N}{N+2}+h_\eps <2\frac{N-1}{N}\leq r_2,
$$
where we use again \eqref{hyp2}. Therefore \eqref{s12.0} is satisfied and then we are done with $\om_2$ as well. At last we discuss the case,
$$
\om_3=\{\frac{1}{r_1}+\frac{1}{r_2}<\frac{N}{N-1},\,r_1\in [r^*_{_N},\frac{2N}{N+2}),\,\,r_2\leq\frac{N}{2}\}.
$$
If $(r_1,r_2)\in\om_3$, since $r_1\in [r^*_{_N},\frac{2N}{N+2})$, then for some $\sg_\eps\to 0$, $h_\eps\to 0$, as $\eps\to 0^+$,
we have that if $\frac{1}{s_1}=\frac{1}{r_1}+\eps$ then $r_1>s_1>2\frac{2N}{N+2}-\sg_\eps$, whence in particular
$$
s_2<2\frac{N-1}{N}+h_\eps <r^{**}_{_N}\leq r_2,
$$
where $r^{**}_{_N}$ is the intersection of $\frac{1}{r_1}+\frac{1}{r_2}<\frac{N}{N-1}$ with the line $r_1=\frac{N}{2}$.
Therefore \eqref{s12.0} is satisfied in this case as well, which concludes the proof of the claim of STEP 1.\\

STEP 2. We prove the existence of the minimum of $J_{\ssl}$ by the direct method.\\
Let $(\rh_{1,n},\rh_{2,n})$ be a minimizing sequence, by STEP 1 $\|\rh_{1,n}\|_{r_1}$ and $\|\rh_{2,n}\|_{r_2}$ are bounded and
we assume that $\rh_{2,n}$ weakly converges in $L^{r_2}(\om)$ and $\rh_{1,n}$ weakly converges in $L^{r_1}(\om)$ to some
$(\rh_{1,\ssl},\rh_{2,\ssl})$. Clearly ${\scriptstyle \frac{1}{r_1}}\ino (\rh_1)^{r_1}+{\scriptstyle \frac{1}{r_2}}\ino (\rh_2)^{r_2}$
is convex whence lowersemicontinuous with respect to the weak topology in $L^{r_1}(\om)\times L^{r_2}(\om)$. We will conclude the proof
showing that, possibly along a subsequence,  we have
$$
\lim\limits_{n\to +\ii}\ino \rh_{1,n}G[\rh_{2,n}]=\ino\rh_{1,\ssl}G[\rh_{2,\ssl}].
$$
Indeed, writing
$$
\ino \rh_{1,n}G[\rh_{2,n}]-\rh_{1,\ssl}G[\rh_{2,\ssl}]=\ino \rh_{1,n}(G[\rh_{2,n}-\rh_{2,\ssl}])+\ino (\rh_{1,n}-\rh_{1,\ssl})G[\rh_{2,\ssl}]=
$$
$$
\ino \rh_{1,n}(G[\rh_{2,n}-\rh_{2,\ssl}])+\ino \rh_{2,\ssl}(G[\rh_{1,n}-\rh_{1,\ssl}]),
$$
it is enough to prove that the embeddings of $W^{2,r_2}(\om)$ in $L^{r^{'}_1}(\om)$ and of $W^{2,r_1}(\om)$ in $L^{r^{'}_2}(\om)$ are compact,
where $r_i^{'}$ is the exponent conjugate to $r_i$. Clearly it is just enough to prove the former.
By \eqref{Souto} we readily deduce that,
$$
r_2>\frac{N-r_1-1}{N(r_1-1)+1},
$$
whence by the Sobolev embedding we find that $W^{2,r_2}(\om)$ is compactly embedded in $W^{1,t}(\om)$ for any
$$
t<t_2=r_1^{'}\frac{N^2-N}{N^2-N+r_1^{'}}=r_1^{'}-\frac{r_1^{'}}{N^2-N+r_1^{'}},
$$
where we remark that $N^2-N+r_1^{'}>r_1^{'}>0$. Thus, again by the Sobolev embedding, we find that
$W^{2,r_2}(\om)$ is compactly embedded in $L^{k}(\om)$ for any
$$
k<k_2=\frac{N r_1^{'}\left(1-\frac{1}{N^2-N+r_1^{'}}\right)}{N-r_1^{'}\left(1-\frac{1}{N^2-N+r_1^{'}}\right)},
$$
where we assume without loss of generality that $N-r_1^{'}\left(1-\frac{1}{N^2-N+r_1^{'}}\right)>0$.
Therefore it is enough to prove that,
$$
\frac{N r_1^{'}\left(1-\frac{1}{N^2-N+r_1^{'}}\right)}{N-r_1^{'}\left(1-\frac{1}{N^2-N+r_1^{'}}\right)}>r_1^{'},
$$
which, after a straightforward evaluation, takes the form,
$$
N^2-\frac{1}{r_1}N+\frac{1}{r_1-1}>0.
$$
The determinant of this polynomial is $\frac{1}{r_1^2}-\frac{4}{r_1-1}$ which is readily seen to be always negative.
Thus we have proved the existence of at least one minimizer of $J_{\ssl}$.

\bigskip

STEP 3 We prove that any minimizer $(\rh_{1,\ssl},\rh_{2,\ssl})$ defines a solution of $\prl$.\\
Let $\om_+=\{x\in\om\,:\rh_{2,\ssl}>0,\mbox{ a.e.}\}$ and $\om_{0}=\{x\in \om\,:\rh_{2,\ssl}=0,\mbox{ a.e.}\}$.
Since $\ino \rh_{2,\ssl}=1$, then $0<|\om_+|\leq |\om|$. For any $n\in \N$ and for any variation of the
form $(\rh_{1,\ssl},\rh_{2,\ssl} +\eps \eta)$,
with supp$(\eta)\subseteq \{x\in\om_+:\rh_{2,\ssl}>\frac1n, \mbox{ a.e.}\}$,
$\eta\in L^{\ii}(\om)$ and $\ino \eta=0$, by using the minimality of $(\rh_{1,\ssl},\rh_{2,\ssl})$ we have,
$$
\int\limits_{\{x\in\om_+:\rh_{2,\ssl}>\frac1n\}}
\left((\rh_{2,\ssl})^\frac{1}{p_2}-\lm G[\rh_{1,\ssl}]\right)\eta\geq o(1), \mbox{ as }\eps\to 0.
$$
Since this is true for any such $\eta$ with $\ino \eta=0$ and any $n\in \N$, then we have that,
\beq\label{eul1}
(\rh_{2,\ssl}(x))^\frac{1}{p_2}-\lm G[\rh_{1,\ssl}](x)=\al_2 \mbox{ a.e. in }\om_+,
\eeq
for a suitable constant $\al_2\in\R$. Next, let $\chi_{A}$ denote the characteristic function of the set $A$,
and assume that $|\om_0|>0$. For any variation of the form $(\rh_{1,\ssl},\rh_{2,\ssl} +\eps \eta)$, with
$$
\eta = \varphi \chi_{\om_{0}}-\left(\,\int\limits_{\om_{0}}\varphi\right)\frac{\chi_{\om_+}}{|\chi_{\om_+}|},
\quad\eta\in L^{\ii}(\om),\;\varphi\geq 0,
$$
we have,
$$
\int\limits_{\om_{0}} \left(-\al_2-\lm G[\rh_{1,\ssl}]\right)\varphi\geq o(1), \mbox{ as }\eps\to 0.
$$
Therefore we conclude that,
\beq\label{eul2}
\al_2+\lm G[\rh_{1,\ssl}](x)\leq 0\mbox{ a.e. in }\om_{0},
\eeq
and in particular that $\psi_{2,\lm}=G[\rh_{1,\ssl}]\in W^{2,r_1}_0(\om)$ is a strong solution of
the first equation in $\prl$ with $\psi_{1,\lm}=G[\rh_{2,\ssl}]\in W^{2,r_2}_0(\om)$, where $(\al_2+\lm\psi_{2,\lm})_+^{p_2}\equiv(\al_2+\lm\psi_{2,\lm})^{p_2}$
as far as $|\om_0|=0$. The same argument shows that $\psi_{1,\lm}$
is a strong solution of the second equation in $\prl$.

\bigskip

At last, by Lemma \ref{lemE1} any strong solution determined in STEP 3 is a classical solution. This fact concludes the proof.
\finedim

\bigskip
\bigskip

\bte\label{monot1} For any $\lm\in [0,\lm^*(\om,{\bf p}))$
$F_{\ssl}$ is real analytic, decreasing with $\frac{d F_{\ssl}}{d\lm}<0$, and concave and we have
$F_{\ssl}^{'}=-\el$. In particular $\frac{d \el}{d\lm}\geq 0$ and $\frac{d}{d\lm}
\left(\frac{p_1\alv[1]}{p_1+1}+\frac{p_2\alv[2]}{p_2+1}\right)< 0$ for any $\lm\in [0,\lm^*(\om,{\bf p}))$.
\ete
\proof
We first observe that, by the uniqueness in Lemma \ref{unique}, we have that $F_{\ssl}$ is the same as $F(\lm)$ in {\bf (VP)}.
It is useful at this stage to introduce the vectorial density,
$$\ur=(\rh_1,\rh_2)$$ and the corresponding entropy,
$$
\mathcal{S}(\ur\,)={\scriptstyle \frac{1}{r_1}}\ino (\rh_1)^{r_1}+
{\scriptstyle \frac{1}{r_2}}\ino (\rh_2)^{r_2},
$$
and energy
$$
\mathcal{E}(\ur\,)=\ino \rh_1 G[\rh_2],
$$
whence ${J}_{\ssl}(\ur\,)=\mathcal{S}(\ur\,)-\lm\mathcal{E}(\ur\,)$.\\
By Lemma \ref{unique} $F_{\ssl}$ is real analytic. If $\lm_2>\lm_1\geq 0$,
then  $J_{\sscp \lm_2}(\ur\,)<J_{\sscp \lm_1}(\ur\,)$, whence
$F_{\sscp \lm_2}\leq F_{\sscp \lm_1}$. Thus
$F_{\ssl}$ is decreasing for $\lm>0$. Moreover, setting
$\lm=t\lm_1+(1-t)\lm_2$, $t\in [0,1]$, and letting $\urls$ be any minimizer of ${J}_{\ssl}$, we find that,
$$
F_{\ssl}=\mathcal{S}(\urls\,)-\lm\mathcal{E}(\urls\,)=
-\mathcal{S}(\urls\,)(t+(1-t))-(t\lm_1+(1-t)\lm_2)\mathcal{E}(\urls\,)=
$$
$$
t J_{\sscp \lm_1}(\urls)+(1-t) J_{\sscp \lm_2}(\urls)\geq t F{\sscp \lm_1}+(1-t)F{\sscp \lm_2}.
$$
Therefore $F_{\ssl}$ is concave with
$\frac{d^2 F_{\ssl}}{d\lm^2}\leq 0$.\\

At this point we use a well known trick about canonical variational principles (see for example \cite{clmp2}).
Let $\lm_1\neq\lm_2$ in $[0, \lm^*(\om,{\bf p}))$ and let $\urlu,\urld$ be the minimizers
of $J_{\sscp \lm_1},J_{\sscp \lm_2}$ respectively. Clearly we have
\beq\label{fb1}
F_{\sscp \lm_1}\leq \mathcal{S}(\urld\,)-\lm_1\mathcal{E}(\urld\,)=
F_{\sscp \lm_2}-(\lm_1-\lm_2)\mathcal{E}(\urld\,),
\eeq
\beq\label{fb2}
F_{\sscp \lm_2}\leq \mathcal{S}(\urlu\,)-\lm_2\mathcal{E}(\urlu\,)=
F_{\sscp \lm_1}-(\lm_2-\lm_1)\mathcal{E}(\urlu\,).
\eeq

and we deduce from \rife{fb1}, \rife{fb2} that

$$
-\mathcal{E}(\urlu\,)\leq \frac{F_{\sscp \lm_1}-F_{\sscp \lm_2}}{\lm_1-\lm_2}\leq
-\mathcal{E}(\urld\,) \quad  \mbox{ if }\lm_1>\lm_2,
$$
$$
-\mathcal{E}(\urld\,)\leq \frac{F_{\sscp \lm_2}-F_{\sscp \lm_1}}{\lm_2-\lm_1}\leq
-\mathcal{E}(\urlu\,)\quad \mbox{ if }\lm_2>\lm_1.
$$

By Lemma \ref{unique}, as $\lm_2\to \lm_1$ we have $\urld \to \urlu$ smoothly, whence

$$
\frac{d F_{\sscp \lm_1}}{d\lm}=-\mathcal{E}(\urlu\,)\equiv- \left.\el\right|_{\lm=\lm_1}.
$$
Remark that $\el>0$ in $[0, \lm^*(\om,{\bf p}))$, whence $\frac{d F_{\ssl}}{d\lm}<0$.
In particular, since $F_{\ssl}$ is real analytic and
concave, then $\frac{d \el}{d\lm}=-\frac{d^2 F_{\ssl}}{d\lm^2}\geq 0$, for any $\lm \in [0, \lm^*(\om,{\bf p}))$.\\
At this point observe that,

$$
F_{\ssl}=\frac{1}{r_1}\ino (\rlv[1])^{r_1}+\frac{1}{r_2}\ino (\rlv[2])^{r_2}-\lm
\ino \rlv[1] G[\rlv[2]]=
$$
$$
 \frac{1}{r_1}\ino\rlv[1](\alv[1]+\lm_1\plv[1])  + \frac{1}{r_2}\ino\rlv[2](\al_2+\lm_2\plv[2])-\lm
\ino \rlv[1] G[\rlv[2]]=
$$
$$
\frac{\alv[1]}{r_1}+\frac{\lm}{r_1}\ino\rlv[1] G[\rlv[2]]+ \frac{\alv[2]}{r_2}+\frac{\lm}{r_2}
\ino\rlv[2] G[\rlv[1]]-\lm
\ino \rlv[1] G[\rlv[2]]=
$$
$$
  \frac{p_1\alv[1]}{p_1+1}+  \frac{p_2\alv[2]}{p_2+1}+\frac{p_1p_2-1}{(p_2+1)(p_1+1)}\el.
$$

By Lemma \ref{unique}, $\alv[i]$, $i=1,2$ are real analytic, and then we deduce that,

$$
\frac{d }{d\lm}\left(\frac{p_1\alv[1]}{p_1+1}+  \frac{p_2\alv[2]}{p_2+1}\right)=
\frac{d F_{\ssl}}{d\lm}-\frac{p_1p_2-1}{(p_2+1)(p_1+1)}\frac{d \el}{d\lm}< 0,
$$
where the strict inequality follows from $\frac{d F_{\ssl}}{d\lm}<0$.\finedim

\bigskip
\bigskip

Finally we prove the estimate about the derivative of $\el$.
\bpr\label{pr-enrg} Let $(\bal,\bpl)\in \mathcal{G}(\om)$ be the unique positive solutions of\, {\rm $\prl$}
for $\lm\in [0, \lm^*(\om,{\mathbf p}))$. Then,
{\rm
\beq\label{8.12.11}
\frac{d \el}{d \lm} \geq p_1\|[\plv[1]]_{1,\ssl}\|^2_{1,\ssl}+p_2\|[\plv[2]]_{2,\ssl}\|^2_{2,\ssl}.
\eeq}
\epr

\proof
By Lemma \ref{unique}
$(\bal,\bpl)$ is a real analytic function of $\lm$ and then, by standard elliptic estimates, we have that
$\bvl\in (C^{2}_0(\ov{\om}\,))^2$ where
$$
\bvl=(\vlv[1],\vlv[2])=\dfrac{d \bpl}{d \lm}=
\left(\dfrac{d\plv[1]}{d \lm},\dfrac{d\plv[2]}{d \lm}\right),
$$
is a classical solution of,
$$
\graf{
-\Delta \vlv[1]=\tlv[2] \rlqv[2] \vlv[2]+p_2\rlqv[2]\plv[2]+p_2 \rlqv[2] \frac{d\alv[2]}{d\lm},\\
-\Delta \vlv[2]=\tlv[1] \rlqv[1] \vlv[1]+p_1\rlqv[1]\plv[1]+p_1 \rlqv[1] \frac{d\alv[1]}{d\lm},
}
$$
and $\frac{d\alv[i]}{d\lm}$ can be computed by the unit mass constraints in $\prl$, that is
$$
p_2\frac{d\alv[2]}{d\lm}=-\tlv[2]  <\vlv[2]>_{2,\ssl}-p_2<\plv[2]>_{2,\ssl},\;
p_1\frac{d\alv[1]}{d\lm}=-\tlv[1]  <\vlv[1]>_{1,\ssl}-p_1<\plv[1]>_{1,\ssl}.\;
$$
Therefore, we conclude that $\bvl$ is a solution of,
\beq\label{1b1}
\graf{
-\Delta \vlv[1]=\tlv[2] \rlqv[2] [\vlv[2]]_{2,\ssl}+p_2\rlqv[2][\plv[2]]_{2,\ssl},\\
-\Delta \vlv[2]=\tlv[1] \rlqv[1] [\vlv[1]]_{1,\ssl}+p_1\rlqv[1][\plv[1]]_{1,\ssl}.
}
\eeq
Since $\el=\ino (\nabla \plv[1],\nabla \plv[2])$,
then by using $\prl$ and \rife{1b1}
we also have,
$$
\frac{d}{d \lm}\el=\ino (\nabla \vlv[1],\nabla \plv[2])+\ino (\nabla \plv[1],\nabla \vlv[2])=
$$
$$
\sum\limits_{i=1}^2\tlv[i]<[\vlv[i]]_{i,\ssl}[\plv[i]]_{i,\ssl}>_{i,\ssl}+
\sum\limits_{i=1}^2p_i\|[\plv[i]]_{i,\ssl}^2\|^2_{i,\ssl}.
$$

Let us denote, for $i=1,2$,
$$
[\plv[i]]_{i,\ssl}=\sum\limits_{k=1}^{+\ii}\xi_{i,k}[\phi_{i,k}]_{i,\ssl}, \quad
[\vlv[i]]_{i,\ssl}=\sum\limits_{k=1}^{+\ii}\beta_{i,k}[\phi_{i,k}]_{i,\ssl},
$$
$$
\xi_{i,k}=<[\phi_{i,k}]_{i,\ssl}[\plv[i]]_{i,\ssl}>_{i,\ssl}, \quad
\beta_{i,k}=<[\phi_{i,k}]_{i,\ssl}[\vlv[i]]_{i,\ssl}>_{i,\ssl},
$$
the Fourier expansions of $[\plv[i]]_{i,\ssl}$ and $[\vlv[i]]_{i,\ssl}$ in $Y_{i,0}$ (see Lemma \ref{ortbase}),
with
respect to the normalized eigenfunctions $[\phi_{i,k}]_{i,\ssl}$, satisfying
$\|[\phi_{i,k}]_{i,\ssl}\|_{i,\ssl}=1$, $i=1,2$. Then, by Lemma \ref{ortbase}, we have,
\beq\label{del}
\frac{d}{d \lm}\el =
\tlv[2]\sum\limits_{k=1}^{+\ii}\beta_{2,k}\xi_{2,k}+
\tlv[1]\sum\limits_{k=1}^{+\ii}\beta_{1,k}\xi_{1,k}+
p_2\sum\limits_{k=1}^{+\ii}\xi^2_{2,k}+p_1\sum\limits_{k=1}^{+\ii}\xi^2_{1,k}.
\eeq
We can now consider $\lambda>0$. On the other side, multiplying the first equation in \rife{1b1} by $\phi_{2,k}$, the second
by $\phi_{1,k}$, using \rife{lineqm} and integrating by parts, we have,

\beq\label{lamq31}
\graf{(\tlv[1]+p_1\sg_k)\beta_{1,k}=\tlv[2]\beta_{2,k}+p_2\xi_{2,k},\\
(\tlv[2]+p_2\sg_k)\beta_{2,k}=\tlv[1]\beta_{1,k}+p_1\xi_{1,k},}
\eeq
where $\sg_{k}=\sg_{k}(\bal,\bpl)$.
Since by \rife{lineq2} we also have $2\lm+\sg_k>\lm$, then \rife{lamq31}
admits the unique solution,
$$
\beta_{1,k}=\frac{\lm(p_1\xi_{1,k}+p_2\xi_{2,k})+\sg_kp_2\xi_{2,k}}{p_1\sg_k(2\lm+\sg_k)},\;
\beta_{2,k}=\frac{\lm(p_1\xi_{1,k}+p_2\xi_{2,k})+\sg_kp_1\xi_{1,k}}{p_2\sg_k(2\lm+\sg_k)},
$$
which we can substitute in \rife{del} to deduce that
\beq\label{enlast}
\frac{d}{d \lm}\el =
\lm\sum\limits_{k=1}^{+\ii}
\frac{\lm(p_1\xi_{1,k}+p_2\xi_{2,k})(\xi_{1,k}+\xi_{2,k})+2\sg_k\xi_{1,k}\xi_{2,k}}{\sg_k(2\lm+\sg_k)}+
p_2\sum\limits_{k=1}^{+\ii}\xi^2_{2,k}+p_1\sum\limits_{k=1}^{+\ii}\xi^2_{1,k}.
\eeq

At this point observe that the equations in $\prl$ can be written in the following form,

$$
\graf{-\Delta \plv[1] =\rlqv[2] (\alv[2]+{\lm}\plv[2])\quad \mbox{in}\;\;\om\\ \\
-\Delta \plv[2] =\rlqv[1] (\alv[1]+{\lm}\plv[1])\quad \mbox{in}\;\;\om\\ \\
\bigintss\limits_{\om} \rlqv[i] (\alv[i]+{\lm}\plv[i])=1\quad i=1,2.
}
$$

Therefore we can evaluate

$$
\alv[i]=\frac{1}{\mlv[i]}-\lm<\plv[i]>_{i,\ssl},\quad i=1,2,
$$

and deduce that

\beq\label{2701.1}
\graf{-\Delta \plv[1] =(\mlv[2])^{-1}\rlqv[2] +{\lm}\rlqv[2][\plv[2]]_{2,\ssl}\quad \mbox{in}\;\;\om\\ \\
-\Delta \plv[2] =(\mlv[1])^{-1}\rlqv[1] +{\lm}\rlqv[1][\plv[1]]_{1,\ssl}\quad \mbox{in}\;\;\om\\ \\
\plv[i]=0 \quad \mbox{on}\;\;\pa\om,\quad i=1,2.
}
\eeq

Multiplying the first equation in \rife{2701.1} by $\phi_{2,k}$, the second
by $\phi_{1,k}$, using \rife{lineqm} and integrating by parts, we have,

\beq\label{lamq31.1}
\graf{(\tlv[1]+p_1\sg_k)\xi_{1,k}=<\phi_{2,k}>_{2,\ssl}+\lm\xi_{2,k},\\
(\tlv[2]+p_2\sg_k)\xi_{2,k}=<\phi_{1,k}>_{1,\ssl}+\lm\xi_{1,k}.}
\eeq
Since $p_i\geq 1$, $i=1,2$ and since by \rife{lineq2} we also have $2\lm+\sg_k>\lm$, then
\rife{lamq31.1} admits the unique solution,
$$
\xi_{1,k}=\frac{(\tlv[2]+p_2\sg_k)<\phi_{2,k}>_{2,\ssl}+\lm<\phi_{1,k}>_{1,\ssl}}{p_1p_2(\lm+\sg_k)^2-\lm^2},\;
$$
$$
\xi_{2,k}=\frac{(\tlv[1]+p_1\sg_k)<\phi_{1,k}>_{1,\ssl}+\lm<\phi_{2,k}>_{2,\ssl}}{p_1p_2(\lm+\sg_k)^2-\lm^2},
$$
and since of course we can assume that $<\phi_{i,k}>_{i,\ssl}\geq 0$, $i=1,2$, for any $k\in\N$,
then we deduce that
\beq\label{enlast1}
\xi_{i,k}\geq 0, i=1,2, \forall\, k\in\N.
\eeq
Remark that all the relations above, including \rife{enlast}, \rife{enlast1}, hold just
by assuming that
$0\notin \sigma({\mathbf L}_{\ssl})$. However at this point we use the fact that for any
$\lm<\lm^*(\om,{\mathbf p})$ it holds $\sg_k\geq \sg_1>0$, $\forall\, k\in\N$, which, together with
\rife{enlast} and \rife{enlast1} implies that

$$
\frac{d}{d \lm}\el \geq
p_2\sum\limits_{k=1}^{+\ii}\xi^2_{2,k}+p_1\sum\limits_{k=1}^{+\ii}\xi^2_{1,k}=
p_2\|[\plv[2]]_{2,\ssl}\|^2_{2,\ssl}+p_1\|[\plv[1]]_{1,\ssl}\|^2_{1,\ssl}.
$$
which is \rife{8.12.11}.
\finedim

\bigskip
\bigskip

At last we present the proofs of Theorems \ref{thmLE} and \ref{thmH}.\\
{\it The proof of Theorem \ref{thmLE}}.
Uniqueness and regularity follow immediately from Lemma \ref{unique} and then
a straightforward evaluation yields the behavior in the claim as $\lm\to 0$.
The inequality about the energy for $\lm=0$ is a well known torsional inequality, see \cite{CRa}.
The monotonicity in the claim follows from Theorem \ref{monot1} and Proposition \ref{pr-enrg}.
Finally, if either $\sg_{1,*}=0$ or if $\al_{i,*}=0$, $i=1,2$, it follows from
Propositions \ref{prsigma} and \ref{prsigma1} that $\lm^*(\om,{\mathbf p})\geq \frac{1}{p_2}\Lambda(\om,2p_2)$.
\finedim

\bigskip

{\it The proof of Theorem \ref{thmH}}.
For $(\bmu_{\ssl}, {\bf u}_{\ssl})$ as defined in the claim, from the constrains in $\prl$ we have
\beq\label{alu}
\alv[i]=\|1+u_{i,\ssl}\|^{-1}_{p_i},\, i=1,2,
\eeq
which immediately shows that
$$
\frac{p_1\alv[1]}{p_1+1}+\frac{p_2\alv[2]}{p_2+1}=\gamma(\bmu,{\bf u}_{\ssl}).
$$

Next observe that
$$
\el=\frac12 \ino \rlv[1]\plv[1]+\frac12 \ino \rlv[2]\plv[2]=
\frac{\alv[1]^{p_1+1}}{2\lm}\ino (1+u_1)^{p_1}u_1+
\frac{\alv[2]^{p_2+1}}{2\lm}\ino (1+u_2)^{p_2}u_2,
$$
which together with \rife{alu} and
\beq\label{lmu}
\lm = \mu_{1,\ssl}\frac{\alv[2]}{\alv[1]^{p_1}},\quad
\lm = \mu_{2,\ssl}\frac{\alv[1]}{\alv[2]^{p_2}},
\eeq
yields
$$
\el=\frac{\alv[1]^{p_1}\alv[2]^{p_2}}{2\mu_{2,\ssl}}\ino (1+u_1)^{p_1}u_1+
\frac{\alv[2]^{p_2}\alv[1]^{p_1}}{2\mu_{1,\ssl}}\ino (1+u_2)^{p_2}u_2=E(\bmu_{\ssl},{\bf u}_{\ssl}).
$$
Moreover, we immediately have $F(\bmu_{\ssl},{\bf u}_{\ssl})=F_{\ssl}$ as well.

Therefore, as far as $\lm<\lm^*(\om,{\bf p})$ the monotonicity properties in the claim follow
immediately from Theorem \ref{thmLE}.\\
Clearly, as far as $\lm<\lm^*(\om,{\bf p})$, by definition we also have $\sg_{1}(\bal,\bpl)>0$ and
$\alv[i]>0$, $i=1,2$, whence by Lemma \ref{unique} we see that $(\bmu_{\ssl}, {\bf u}_{\ssl})$ is a continuous
real analytic curve.\\
Finally, assume that $\al_{1,*}=0$ and $\al_{2,*}>0$. By definition there exists a
sequence $\lm_n\to (\lm^*(\om,{\mathbf p}))^{-}$, such
that $\al_{1,\sscp \lm_n}\to 0^+$, $\al_{2,\sscp \lm_n}\to \al_2>0$. The curve is obviously unbounded since by \rife{lmu} we have
$\mu_{2,\ssl}\to +\ii$.

Now, under the assumption of Theorem \ref{thlambda}-$(b)$, we have $\lm^*(\om,{\mathbf p})<+\ii$ and then by Lemma \ref{lemE1}
and passing to a further subsequence if necessary we can assume that $\psi_{i,\sscp \lm_n}$, $i=1,2$
converge smoothly to $\psi_{i,*}$, $i=1,2$, which are classical solutions of $\prl$ for
$\lm=\lm^*(\om,{\mathbf p})$ and $\al_1=0$, $\al_2>0$. Since by definition
$\al_{i, \sscp \lm_n}u_{i,\sscp \lm_n}=\lm_n\psi_{i,\sscp \lm_n}$, $i=1,2$, then the convergence in the claim
follows from \rife{alu}.\\
The conclusion in case $\al_{2,*}=0$ and $\al_{1,*}>0$ follows in the same way.
\finedim

\bigskip
\bigskip

\appendix

\section{Uniqueness of solutions for $\lm$ small}\label{appF}
The following lemma is proved by a standard application of the contraction mapping principle and we
prove it here just for reader's convenience. We will denote by $C_1$ the constant in Lemma \ref{lemE1}
and by $C_2,C_3$ other positive constants depending only by $r_0,\om,{\bf p},N$.
\ble\label{lmsmall}
There exists $\lm_0>0$ such that:\\
$(j)$ for any $\lm\in [0,\lm_0]$ there exist at least one solution
$(\bal,\bpl)$ of {\rm $\prl$}.\\
$(jj)$ for any solution of {\rm $\prl$} we have
$\alv[i]> \frac13$, $i=1,2$ for any $\lm\in [0,\lm_0]$.
\ele
\proof $(j)$\\
Putting $u_1=\lm\psi_1,u_2=\lm \psi_2$ the proof is an immediate consequence of the following lemmas.
Let $\ba=(\al_1,\al_2)$, ${\bf u}=(u_1,u_2)$ and let us define,
$$
{\bf B}_{\ii}=\left\{{\bf u}\in (L^{\ii}(\om))^2\,
|\,\|{\bf u}\|_{L^{\ii}(\om)}:=\max\limits_{i=1,2}\|u_i\|_{L^{\ii}(\om)}\leq C_1,\,
u_i\geq 0,\mbox{ a.e. in } \om,\,i=1,2\right\}
$$
where $C_1(r,\om,1,{\bf p},N)$ is the constant obtained in Lemma \ref{lemE1} evaluated with $\ov{\lm}=1$.

\ble\label{lemE2} Let $(p_1,p_2)$ satisfy \eqref{Souto}. There exists $\lm_0\in (0, 1]$
such that for any $\lm\in [0,\lm_0]$ and for any $\al_i\in (-\ii,1]$ there exists
one and only one solution ${\bf u}=(u_{1,\sscp \lm,\ba},u_{2,\sscp \lm,\ba})\in (C^{2,r}_0(\ov{\om}\,))^2$
of the problem
\beq\label{B.1}
\graf{-\Delta u_1 =\lm (\al_2+u_2)_+^{p_2}\quad \mbox{in}\;\;\om\\ \\
-\Delta u_2 =\lm (\al_1+u_1)_+^{p_1}\quad \mbox{in}\;\;\om\\ \\
u_1=0=u_2 \quad \mbox{on}\;\;\pa\om\\\\
u_i\in B_{\ii},\,i=1,2.}
\eeq
Moreover, for fixed $\lm\in [0,\lm_0]$, the maps
$(-\ii,1]\ni \al_1\to u_{2,\lm}[\al_1]=u_{2,\lm,\ba}\in B_{\ii}$,
$(-\ii,1]\ni \al_2\to u_{1,\lm}[\al_2]=u_{1,\lm,\ba}\in B_{\ii}$,
are continuous and $u_{1,\lm,\ba}\equiv 0\equiv u_{2,\lm,\ba}$ if either $\lm=0$ or if $\al_i\leq 0$, $i=1,2$.
\ele
\proof
First of all, if either $\lm=0$ or if $\al_i\leq 0$, $i=1,2$, then  $(u_1,u_2)\equiv (0,0)$ is a solution, whence the last part
of the statement will follow immediately from the uniqueness.\\
For ${\lm}_0\in (0,1]$ to be fixed later on and for fixed $\lm\in [0,\lm_0]$ and $\al_i\in (-\ii,1]$, $i=1,2$,
we define
$$
{\bf T}_{\lm,\ba}({\bf u})=\lm (G[(\al_2+u_2)_+^{p_2}],G[(\al_1+u_1)_+^{p_1}] ),\quad u_i\in B_{\ii}, i=1,2.
$$

Recall that if $\al_i>0$ then $\al_i\leq 1$, while if $\al_i<0$ then $(\al_i+u_i)_+\leq (u_i)_+$,
whence we have,
$$
\|{\bf T}_{\lm,\ba}({\bf u})\|_{L^{\ii}(\om)}:=
\lm\max\limits_{i=1,2} \|G[(\al_i+u_i)_+^{p_i}]\|_{L^{\ii}(\om)}\leq \lm C_2,
$$
and we readily see that
${\bf T}_{\lm,\ba}:{\bf B}_{\ii}\to {\bf B}_{\ii}$ for any $\lm\leq \frac{C_1}{C_2}$. Also,
$$
\|{\bf T}_{\lm,\ba}({\bf u})-{\bf T}_{\lm,\ba}({\bf v})\|_{L^{\ii}(\om)}\leq \lm
\max \limits_{i=1,2}\|p_iG[(\al_i+w_i)_+^{p_i-1}|u_i-v_i|]\|_{L^{\ii}(\om)}\leq
$$
$$
\lm \max \limits_{i=1,2} \|p_iG[(\al_i+w_i)_+^{p_i-1}\|_{L^{\ii}(\om)} \|u_i-v_i\|_{L^{\ii}(\om)}\leq
\lm C_3\|{\bf u}-{\bf v}\|_{L^{\ii}(\om)}
$$
where $w_i\in B_{\ii}$ satisfies $u_i\leq w_i\leq v_i$, $i=1,2$.

Therefore, we also have $\|{\bf T}_{\lm,\ba}({\bf u})-{\bf T}_{\lm,\ba}({\bf v})\|_{L^{\ii}(\om)}\leq
\frac12\|{\bf u}-{\bf v}\|_{L^{\ii}(\om)}$, for any
$\lm\leq \frac{1}{2C_3}$. As a consequence putting $\lm_0= \min\{1,\frac{C_1}{C_2},\frac{1}{2C_3}\}$,
we have that ${\bf T}_{\lm,\ba}$ is a contraction in ${\bf B}_{\ii}$ for any $\lm\leq \lm_0$.
Whence, in particular, for any fixed $\ba\in ((-\ii,1])^2$, we have that for any $\lm\in [0,\lm_0\,]$
there exists a unique solution of ${\bf u}={\bf T}_{\lm,\ba}({\bf u})$. The existence and uniqueness claim follows since,
by standard elliptic estimates, $(u_{1,\sscp \lm,\ba},u_{2,\sscp \lm,\ba})\in (C^{2,r}_0(\ov{\om}))^2$ solves the problem
in the statement of the lemma if and only if ${\bf u}\in {\bf B}_{\ii}$ satisfies ${\bf u}={\bf T}_{\lm,\ba}({\bf u})$.\\
Concerning the continuity of $(u_{1,\lm}[\al_2],u_{2,\lm}[\al_1])=(u_{1,\lm,\ba},u_{2,\lm,\ba})$
for $\al_i\in (-\ii,1],\,i=1,2$, we observe that if $\ba_n=(\al_{1,n},\al_{1,n})\to \ba=(\al_1,\al_2)$, then
$$
\|u_{2,\lm}[\al_{1,n}]-u_{2,\lm}[\al_1]\|_{L^{\ii}(\om)}=
\|G[(\al_{1,n}+u_{1,\lm,\ba_n})_+^{p_1}]-G[(\al_1+u_{1,\lm,\ba})_+^{p_1}]\|_{L^{\ii}(\om)}\leq
$$
$$
\|G[(\al_{1,n}+u_{1,\lm,\ba_n})_+^{p_1}]-G[(\al_{1,n}+u_{1,\lm,\ba})_+^{p_1}]\|_{L^{\ii}(\om)}+
$$
$$
\|G[(\al_{1,n}+u_{1,\lm,\ba})_+^{p_1}]-G[(\al_1+u_{1,\lm,\ba})_+^{p_1}]\|_{L^{\ii}(\om)}\leq
$$
$$
\lm C_3\|u_{1,\lm}[\al_{2,n}]-u_{1,\lm}[\al_2]\|_{L^{\ii}(\om)} +
p_1\lm\|G[(s+u_{1,\lm,\al})_+^{p_1-1}]\|_{L^{\ii}(\om)}|\al_{1,n}-\al_1|\leq
$$

$$
\frac{1}{2}\|u_{1,\lm}[\al_{2,n}]-u_{1,\lm}[\al_2]\|_{L^{\ii}(\om)}+\lm_0 C_3|\al_{1,n}-\al_1|.
$$

Clearly the same argument applies to the other component and then we deduce that
$$
\|u_{2,\lm}[\al_{1,n}]-u_{2,\lm}[\al_1]\|_{L^{\ii}(\om)}+
\|u_{1,\lm}[\al_{2,n}]-u_{1,\lm}[\al_2]\|_{L^{\ii}(\om)}\leq
4\lm_0 C_3(|\al_{1,n}-\al_1|+|\al_{2,n}-\al_2|),
$$
which readily implies the claim.
\finedim

\bigskip
\bigskip

For fixed $\lm\in [0,\lm_0]$ we consider the continuous map
$$
((-\ii,1])^2\ni \ba=(\al_1,\al_2)\to {\bf u}_{\ba} = (u_{1,\lm}[\al_2],u_{2,\lm}[\al_1])\in (B_{\ii})^2,
$$
where $u_{1,\lm}[\al_2]=u_{1,\lm,\ba}, u_{2,\lm}[\al_1]=u_{2,\lm,\ba}$.
Then we have,

\ble\label{lemE3} By taking a smaller $\lm_0$ if necessary, for any fixed $\lm\in [0,\lm_0]$ we have:\\
$(i)$ The maps $u_{1,\lm}[\al_2],u_{2,\lm}[\al_1]$ are monotonic increasing,
$$
u_{1,\lm}[\al_2]\leq u_{1,\lm}[\beta_2],\quad \forall\,0<\al_2<\beta_2\leq 1,\quad
u_{2,\lm}[\al_1]\leq u_{2,\lm}[\beta_1],\quad \forall\,0<\al_1<\beta_1\leq 1.
$$
$(ii)$ There exists at least one  $\bal=(\alv[1],\alv[2])\in ((\frac13,1])^2$ such that,
$$
\ino (\alv[1]+u_{1,\lm}[\alv[2]])^{p_1}=1=\ino (\alv[2]+u_{2,\lm}[\alv[1]])^{p_2}.
$$
\ele
\proof
$(i)$ If $\lm=0$ we have $u_{1,0,\ba}=0$ for any $\ba$ and the conclusion is trivial. For any fixed
$-\ii<\al_2<\beta_2\leq 1$ let us set,
$$
(w_1,w_2)=(u_{1,\lm}[\beta_2]-u_{1,\lm}[\al_2],u_{2,\lm}[\beta_1]-u_{2,\lm}[\al_1])\in
(C^{2,r}_0(\ov{\om}\,))^2,
$$
then
$$
-\Delta w_1=\lm (\beta_2+u_{2,\lm}[\beta_1])_+^{p_2}-\lm(\al_2+u_{2,\lm}[\al_1])_+^{p_2}\geq
\lm (\al_2+u_{2,\lm}[\beta_1])_+^{p_2}-\lm(\al_2+u_{2,\lm}[\al_1])_+^{p_2}\geq
$$
$$
\lm p_2 (\al_2+u_{2,\lm}[\al_1])_+^{p_2-1}(u_{2,\lm}[\beta_1]-u_{2,\lm}[\al_1])=
\lm p_2 (\al_2+u_{2,\lm}[\al_1])_+^{p_2-1}w_2,
$$
by the convexity of $f(t)=(\al+t)_+^p$ for $t\in \R$. By applying the same argument to $w_2$ we deduce that
$$
\graf{ -\Delta w_1 \geq  V_2w_2\\
-\Delta w_2 \geq V_1w_1},\quad V_1=\lm p_1(\al_1+u_{1,\lm}[\al_2])_+^{p_1-1},
\,V_2=\lm p_2(\al_2+u_{2,\lm}[\al_1])_+^{p_2-1}.
$$
and then, in view of Lemma \ref{lemE1} and possibly taking a smaller $\lm_0$,
well known results for cooperative elliptic systems (\cite{dFM}) show that $w_i\geq 0$, $i=1,2$, as claimed.\\

 $(ii)$ For $\lm=0$ we have $u_{i,0,\ba}=0$, $i=1,2$ and then necessarily $\alv[i]=1$, $i=1,2$.
 For fixed $\lm \in (0,\lm_0]$, by Lemma \ref{lemE2} and $(i)$ the functions
$$
g_1(\ba)=\ino (\al_1+u_{1,\lm}[\al_2])_+^{p_1},\quad
g_2(\ba)=\ino (\al_2+u_{2,\lm}[\al_1])_+^{p_2},\quad \ba \in ((-\ii,1])^2,
$$
are continuous as a function of $(\al_1,\al_2)$ and increasing in $\al_1$ and $\al_2$. Clearly
$\|u_{1,\lm}[\al_2]\|_{\ii}\leq \lm C_2$, $\|u_{1,\lm}[\al_2]\|_{\ii}\leq \lm C_2$ for any $\lm\leq \lm_0$, and
then, possibly taking $\lm_0$ small enough to guarantee that
\beq\label{B.0}
(2\lm_0 C_2)^{p_i}\leq \frac{1}{4},\; i=1,2,
\eeq
we have,
$$
g_1((0,\al_2))=
\ino (u_{1,\lm}[\al_2])^{p_1}\leq \frac14,\;\forall \lm \in (0,\lm_0],\;\forall \al_2\in (-\ii,1],
$$
$$
g_2((\al_1,0))\leq
\ino (u_{2,\lm}[\al_1])^{p_2}\leq \frac14,\;\forall \lm \in (0,\lm_0],\;\forall \al_1\in (-\ii,1],
$$

while we also have,
$$
g_1((1,\al_2))=\ino (1+u_{1,\lm}[\al_2])^{p_1} >1,\;\forall\lm \in (0,\lm_0],\;\forall \al_2\in (-\ii,1],
$$
$$
\quad g_2((\al_1,1))=\ino (1+u_{2,\lm}[\al_1])^{p_2} >1,\;\forall \lm \in (0,\lm_0],\;\forall \al_1\in (-\ii,1].
$$
As a consequence, for any $\lm \in (0,\lm_0]$ there exists at least one $\bal=(\alv[1],\alv[2])$ such that
$g_i(\bal)=1$, $i=1,2$ which, by the monotonicity of $g_i$, must necessarily satisfy
\beq\label{B.2}
\alv[i]\in (0,1),\; i=1,2.
\eeq
This concludes the proof of Lemma \ref{lmsmall}-$(j)$.
\medskip

Proof of $(jj)$\\
Let $(\bal,\bpl)$ be any solution of $\prl$ for $\lm\in [0,\lm_0]$. If $\lm=0$ then necessarily $\bal=(1,1)$.
Otherwise let $\lm\in (0,\lm_0]$ and observe that $\lm\bpl={\bf u}_{\ssl}$ where ${\bf u}_{\ssl}$ is the unique
solution of \rife{B.1} found in $(j)$. As a consequence, by \rife{B.0} and \rife{B.2} we have that,
$$
1 =\ino (\alv[1]+u_{1,\lm}[\alv[2]])_+^{p_1}\leq 2^{p_1}(\alv[1])_+^{p_1} +(2 \lm C_2 \,)^{p_1}\leq 2^{p_1}\alv[1]^{p_1}+\frac{1}{4},
$$
whence $\alv[1]\geq (\frac{3}{4})^{\frac1{p_1}}\frac12>\frac13$.
\finedim

\end{document}